\documentclass[12pt]{iopart}
\usepackage[utf8]{inputenc}

\usepackage{graphicx}
\usepackage{iopams}
\usepackage{amsopn}
\usepackage{xcolor}
\usepackage{hyperref}
\usepackage{microtype}
\usepackage{xspace}

\hypersetup{colorlinks=true,allcolors=blue}

\usepackage{textcomp}
\newcommand{\textapprox}{\raisebox{0.5ex}{\texttildelow}}
\newcommand{\iint}{\int \hspace{-3mm} \int}
\newcommand{\iiint}{\int \hspace{-3mm} \int \hspace{-3mm} \int}

\newif\ifdraft

\newcommand{\eqref}[1]{\eref{#1}}
\newcommand{\text}[1]{\mbox{#1}}

\drafttrue

\ifdraft
    \newcommand{\ja}[1]{\textcolor{purple}{\textbf{Joakim: #1}}}
    \newcommand{\jatoar}[1]{\textcolor{purple}{\textbf{Joakim $\rightarrow$ Aaditya: #1}}}
    \newcommand{\jatoms}[1]{\textcolor{purple}{\textbf{Joakim $\rightarrow$ Marina: #1}}}
    \newcommand{\jatoab}[1]{\textcolor{purple}{\textbf{Joakim $\rightarrow$ Alex: #1}}}

    \newcommand{\ar}[1]{\textcolor{cyan}{\textbf{Aaditya: #1}}}
    \newcommand{\artoms}[1]{\textcolor{cyan}{\textbf{Aaditya $\rightarrow$ Marina: #1}}}
    \newcommand{\artoja}[1]{\textcolor{cyan}{\textbf{Aaditya $\rightarrow$ Joakim: #1}}}
    \newcommand{\artoab}[1]{\textcolor{cyan}{\textbf{Aaditya $\rightarrow$ Alex: #1}}}

    \renewcommand{\ms}[1]{\textcolor{blue}{\textbf{Marina: #1}}}
    \newcommand{\mstoja}[1]{\textcolor{blue}{\textbf{Marina $\rightarrow$ Joakim: #1}}}
    \newcommand{\mstoar}[1]{\textcolor{blue}{\textbf{Marina $\rightarrow$ Aaditya: #1}}}
    \newcommand{\mstoab}[1]{\textcolor{blue}{\textbf{Marina $\rightarrow$ Alex: #1}}}

    \newcommand{\ab}[1]{\textcolor{olive}{\textbf{Alex: #1}}}
    \newcommand{\abtoja}[1]{\textcolor{olive}{\textbf{Alex $\rightarrow$ Joakim: #1}}}
    \newcommand{\abtoar}[1]{\textcolor{olive}{\textbf{Alex $\rightarrow$ Aaditya: #1}}}
    \newcommand{\abtoms}[1]{\textcolor{olive}{\textbf{Alex $\rightarrow$ Marina: #1}}}
\else
    \newcommand{\ja}[1]{}
    \newcommand{\jatoar}[1]{}
    \newcommand{\jatoms}[1]{}
    \newcommand{\jatoab}[1]{}

    \newcommand{\ar}[1]{}
    \newcommand{\artoms}[1]{}
    \newcommand{\artoja}[1]{}
    \newcommand{\artoab}[1]{}

    \renewcommand{\ms}[1]{}
    \newcommand{\mstoja}[1]{}
    \newcommand{\mstoar}[1]{}
    \newcommand{\mstoab}[1]{}

    \newcommand{\ab}[1]{}
    \newcommand{\abtoja}[1]{}
    \newcommand{\abtoar}[1]{}
    \newcommand{\abtoms}[1]{}
\fi

\newcommand{\imunit}{\mathrm{i}}
\newcommand{\euler}{\mathrm{e}}
\newcommand{\Real}{\mathbb{R}}
\newcommand{\Complex}{\mathbb{C}}
\newcommand{\Integer}{\mathbb{Z}}

\newcommand{\Leb}{L}

\newcommand{\hs}{\mathrm{HS}}

\newcommand{\vect}[1]{\boldsymbol{#1}}
\newcommand{\fourier}[1]{\hat{#1}}
\newcommand{\bessel}[1]{\lowercase{#1}}
\newcommand{\wfourier}[1]{\widehat{#1}}

\newcommand{\vx}{{\vect{x}}}
\newcommand{\vk}{{\vect{k}}}
\newcommand{\vd}{{\vect{\delta}}}

\newcommand{\fA}{{\fourier{A}}}
\newcommand{\fB}{{\fourier{B}}}

\newcommand{\bA}{{\bessel{A}}}
\newcommand{\bB}{{\bessel{B}}}

\newcommand{\kmax}{{K}}
\newcommand{\dmax}{{D}}
\newcommand{\qmax}{{Q}}
\newcommand{\ellmax}{{L}}
\newcommand{\mmax}{{M}}
\newcommand{\wmax}{{W}}
\newcommand{\dx}{{\Delta x}}     

\newcommand{\argmin}{\operatorname*{arg\,min}}

\newcommand{\bigO}{\mathcal{O}}

\newcommand{\lin}{{\mathrm{lin}}}
\newcommand{\ls}{{\mathrm{LS}}}
\newcommand{\sls}{{\mathrm{GLS}}}
\newcommand{\svd}{{\mathrm{SVD}}}

\newcommand{\ntrans}{N}      
\newcommand{\nnodes}{{H}}
\newcommand{\nnodessls}{{H}}
\newcommand{\nnodesell}{{H_{\ell}}}
\newcommand{\nim}{N_I}               

\newcommand{\besselker}{\mathcal{M}}

\newcommand{\meth}{FTK\xspace}        

\newtheorem{remark}{Remark}
\newtheorem{thm}{Theorem}
\newtheorem{lem}{Lemma}

\begin{document}

\title[FTK for fast rigid image alignment]{Factorization of the translation kernel for fast rigid image alignment}
\author{Aaditya Rangan,$^{1}$ Marina Spivak,$^2$ Joakim And\'en,$^2$ and Alex Barnett$^2$}

\address{$^1$ Courant Institute, New York University, New York, NY and Center for Computational Mathematics, Flatiron Institute, New York, NY} \ead{rangan@cims.nyu.edu}

\address{$^2$ Center for Computational Mathematics, Flatiron Institute, New York, NY} \eads{\mailto{mspivak@flatironinstitute.org}, \mailto{janden@flatironinstitute.org}, \mailto{abarnett@flatironinstitute.org}}

\begin{abstract}
An important component of many image alignment methods is the calculation of inner products (correlations) between an image of $n\times n$ pixels and
another image translated by some shift and rotated by some angle.
For robust alignment of an image pair, the number of considered shifts and angles is typically high, thus the inner product calculation becomes a bottleneck.
Existing methods, based on fast Fourier transforms (FFTs),
compute all such inner products with computational complexity $\bigO(n^3 \log n)$ per image pair, which is reduced to $\bigO(\ntrans n^2)$ if only $\ntrans$ distinct shifts are needed.
We propose to use a factorization of the translation kernel (\meth), an optimal interpolation method which represents images in a Fourier--Bessel basis and uses a rank-$H$ approximation of the translation kernel via an operator singular value decomposition (SVD).
Its complexity is $\bigO(Hn(n + \ntrans))$ per image pair.
We prove that $H = \bigO((W + \log(1/\epsilon))^2)$,
where $2W$ is the magnitude of the maximum desired shift in pixels and $\epsilon$ is the desired accuracy.
For fixed $W$ this leads to an acceleration when $\ntrans$ is large, such as when sub-pixel shift grids are considered.
Finally, we present numerical results in an electron cryomicroscopy application showing speedup factors of $3$--$10$ with respect to the state of the art.
\end{abstract}

\vspace{2pc}

\submitto{\IP}

\maketitle

\section{Introduction}
\label{sec:intro}

Rigid alignment of images is a ubiquitous task that arises in many computer vision problems, such as motion tracking, video stabilization, summarization, image stitching and the creation of mosaics \cite{Szeliski2006}.
It is also widely used in the analysis of biomedical data.
One such example, which serves as the motivation for this paper, is single particle electron cryomicroscopy (cryo-EM). For an overview of cryo-EM, see \cite{Cheng2015,Murata2017,Sigworth2016,Nogales2015,Elmlund2015}.

Within this context, the basic task is to compare (i.e., correlate) pairs of 2D images: each image is sampled on a uniform square pixel grid, and the goal is to find the rigid transformation---a rotation followed by a translation---of one image which maximizes its consistency with the other.
This is necessary at several stages in the cryo-EM pipeline \cite{Wang2013,Singer2009,Shkolnisky2012,Goncharov1987,vanHeel1987,Jonic2005,Yang2008,Zhao2014,Lederman2019}, and plays a major role in many popular software packages \cite{Punjani2017,Scheres2012,Grigorieff2016,Lyumkis2013,Grigorieff2007,EMAN2}.
Outside of cryo-EM, the rigid alignment problem also occurs in certain extensions of non-local means denoising \cite{nonlocalmeans1,nonlocalmeans2} where patches are aligned rotationally instead of just translationally \cite{rotnonlocalmeans1,rotnonlocalmeans2}.

\begin{figure}
\centering
  \includegraphics[width=.9\linewidth,trim= {0cm 0cm 0cm 0cm},clip]{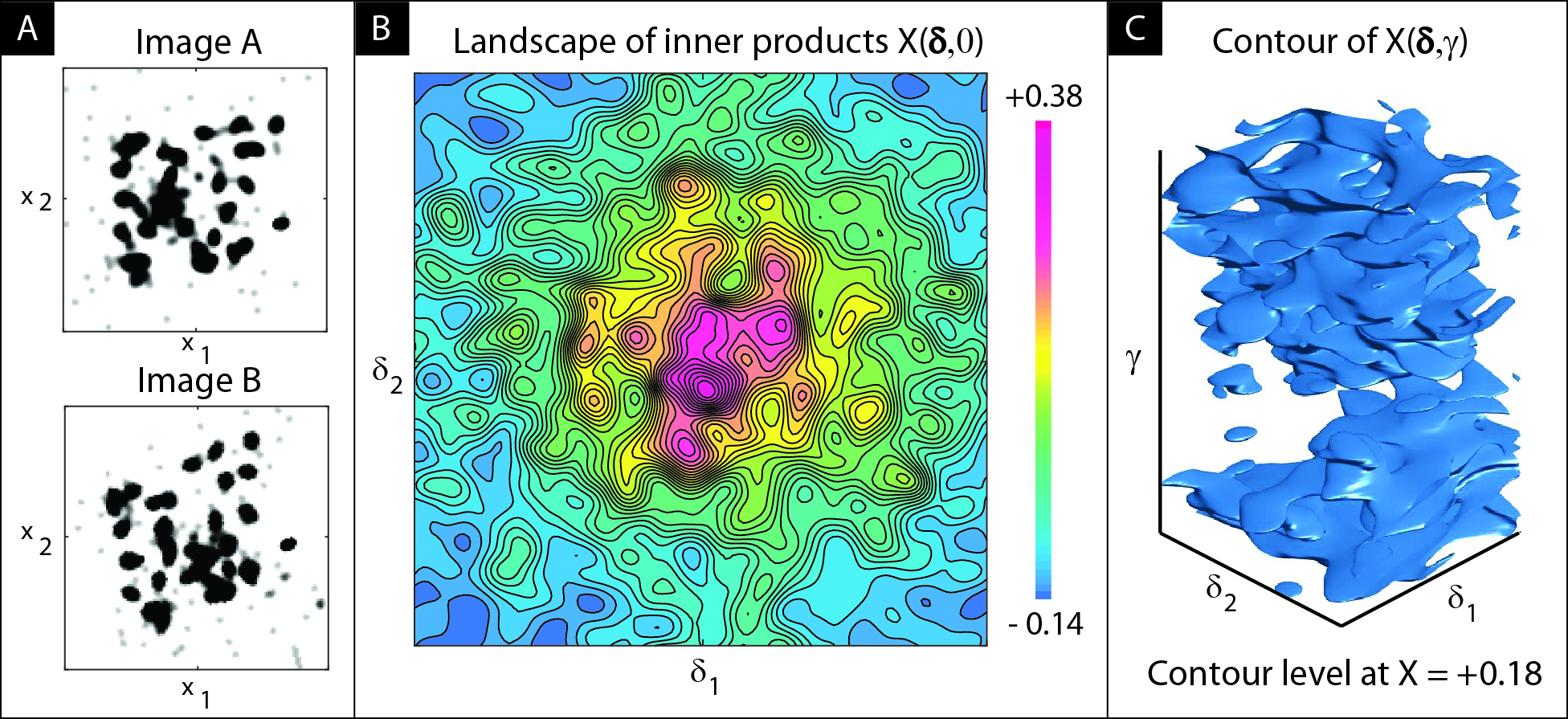}
\caption{\label{fig:landscape}
In panel A we show two images, $A(\vx)$ and $B(\vx)$, for $\vx \in [-1,1]^2$.
Each of these images is generated by summing several anisotropic Gaussian functions, and then normalizing so that the resulting image has $\Leb^2$-norm 1.
In panel B we show the inner product function ${\cal X}(\vd,0)$
(see \eqref{Xdg})
obtained by translating the images with respect to one another.
Note that the landscape is quite complicated, with multiple local maxima (see colorbar to the right).
The 3D landscape obtained by including rotations is even more complex: one level set of ${\cal X}(\vd,\gamma)$ is shown in panel C.
}
\end{figure}

We consider the very common setting where the metric of comparison between
images is their {\em inner product}, equivalent (up to normalization) to their
cross-correlation.
The challenge is then to solve the image-alignment problem
to high (often sub-pixel) accuracy, robustly and rapidly.
As Figure~\ref{fig:landscape} illustrates, the inner product, viewed as
a function of shift and angle, may have a complicated landscape with
many local maxima.
This is particularly true in the high-noise setting relevant to
cryo-EM.
Accurate and robust alignment thus
requires computing inner products for a large set of transformations, which is expensive. In the literature, there are two general approaches for tackling this difficulty.
Some methods conduct an exhaustive search over a grid of rotations and translations, calculating inner products for every transformation on this grid.
Others take shortcuts by using various optimization methods to resolve only a portion of the full grid, yielding a higher resolution than permitted by an exhaustive calculation.

In this paper, we present a new strategy for solving this problem,
which combines many of the advantages of the methods discussed above,
and is an optimal interpolation of the inner product as a function
of translation.
We first review some of the correlation-based methods currently in use, and then describe our approach in greater detail.

\paragraph{Exhaustive methods.}
The simplest approach to the 2D image alignment problem is an exhaustive calculation of inner products over a dense grid of shifts and angles.
The high grid density makes direct calculation computationally prohibitive
(see Table~\ref{table:complexity}), and several accelerations have been proposed \cite{Joyeux2002}.

Two widely used accelerations are based on the convolution theorem, which states that cross-correlation in real space is equivalent to multiplication in Fourier space.
The first, which we refer to as ``brute-force rotations'' (BFR), accelerates the correlation across a full image-sized grid of translations through 2D FFT-accelerated convolutions, then steps through all rotated versions of each image via brute force.
The second method, ``brute-force translations'' (BFT), uses the convolution theorem to treat rotations efficiently via 1D FFTs along rotation angles, while stepping through all required translations by brute force \cite{Joyeux2002,Barnett2017}. Their complexities are listed in the middle two rows of Table~\ref{table:complexity}.

All such exhaustive methods have the advantage of resolving the full landscape of inner products on a pre-defined roto-translation grid, ensuring robust maximization in the case of multiple local maxima.
Moreover, sampling this entire landscape is necessary in a Bayesian framework, where an integral of a probability distribution over shifts and angles is used to marginalize the likelihood \cite{Scheres2012,Punjani2017}.
Unfortunately, these methods all become very costly when
sub-pixel accuracy is desired for shifts and angles.

\paragraph{Optimization-based methods.}
In contrast to the exhaustive methods above, other methods try to find the best alignment quickly by iteratively maximizing the inner product
as a function of shift and angle.
Since this does not explore the entire landscape of inner products, it can be much faster.

One strategy alternates between maximization over rotation angles
(fixing the shift) and shifts (fixing the angle) \cite{Penczek1992,Harauz1988}.
Another method conducts an exhaustive search over a coarse grid to find an approximate maximum that is then refined \cite{Yang2008}.

In favorable situations (e.g., high signal-to-noise ratio images with
dominant low-frequency content), these optimization-based methods quickly converge to the best alignment without being tethered to a pre-specified grid.
Unfortunately, these methods may converge to a locally optimal alignment far from the true global optimum.
In addition, because they only sample part of the inner product landscape, these methods are not immediately useful in a Bayesian setting.

\paragraph{Our approach.}
In this paper, we present the ``factorization of the translation kernel'' (\meth) method, a new algorithm which efficiently evaluates inner products
over all rotation angles and over all shifts lying in a predetermined domain.
It also enables exhaustive high density (sub-pixel) sampling at very little extra computational cost.

The essential ingredients in our approach are (i) the convolution theorem applied to rotation angles (as leveraged in BFT) and (ii) a truncated functional SVD of the translation kernel in the \emph{Fourier--Bessel basis}
compatible with this angular convolution.
This allows us to form, and then evaluate,
an optimal low-rank interpolant for the inner product landscape at all rotation angles and all shifts up to a certain magnitude.

The Fourier--Bessel basis we use to represent our images, along with many similar analytical tools, have also been used by Park {\it et al.} \cite{Park2010} for deblurring the class averages of previously aligned 2D images.
This representation was also used by Barnett {\it et al.} \cite{Barnett2017} to assign viewing angles to images within the cryo-EM reconstruction pipeline.

Unlike many other exhaustive algorithms, our algorithm does not demand a pre-specified grid of translations.
Instead, the speed of our algorithm depends on two parameters: the desired accuracy of the calculation and the maximum magnitude of the desired shifts.
As the desired accuracy increases, the rank of the truncated SVD grows, giving a mild increase in overall computational cost.
As the maximum shift magnitude increases, the structure of the corresponding set of translation operators becomes richer, and more terms in the SVD are required to achieve the same level of accuracy, also increasing the computational cost.

Importantly, given a fixed level of accuracy and range of translations, the computational cost of \meth scales with the number of non-negligible terms in the SVD.
Once these terms are calculated and integrated against the given images, they can be used to recover the inner product values for a fine
sub-pixel translation and rotation grid cheaply.
Table~\ref{table:complexity} compares the complexity of
our proposal to the other exhaustive direct, BFR, and BFT methods.

We demonstrate the utility of \meth by applying it
to what is often the bottleneck in cryo-EM reconstruction: the task known as
``template matching.''
This requires the inner product evaluation between all pairs
of single-particle cryo-EM images and templates, over a range of sub-pixel
translations and all in-plane rotations.
We achieve a substantial speedup of a factor $3$--$10$
over BFT or BFR methods for a wide range of translation grid sizes,
while recovering at least two digits of accuracy (sufficient for most cryo-EM
applications).
We also show that, for higher accuracies, \meth
is several orders of magnitude more efficient
than the commonly-used method of
linear interpolation from a translation grid.

The rest of this paper is organized as follows.
Section~\ref{sec:setup} describes basic aspects of Fourier representation
of images, and defines the inner product function.
Section~\ref{sec:bessel} reformulates this in the Fourier--Bessel basis
and explains its numerical discretization.
The main SVD algorithm is explained in the context of optimal interpolation
schemes in Section~\ref{sec:fast}.
In Section~\ref{sec:erroranalysis} we
prove a theorem bounding the rank growth in terms of the maximum translation
in pixels and the desired relative mean squared error in the inner product calculation.
Section~\ref{sec:results} presents numerical results and a comparison with the BFR and BFT methods in the cryo-EM setting.
In Section~\ref{sec:ext} we describe how \meth extends to other 2D kernels and to 3D volume alignment.
We conclude in Section~\ref{sec:conc}.
A Python implementation is available at \url{https://github.com/flatironinstitute/ftk}.

\begin{table}   
  \centering
  \begin{tabular}{l|l}
Algorithm & Effort per image--template pair\\
    \hline
Direct method & $\bigO(\ntrans n^3)$ \\
Brute force rotations (BFR) &  $\bigO(n^3\log n)$ \\
Brute force translations (BFT) &  $\bigO(\ntrans n^2)$ \\
Proposed method (\meth) &  $\bigO(\nnodes n (n + \ntrans))$\\
\hline
\end{tabular}
  \caption{Complexity of some algorithms for rigid image alignment
    by computation of sets of translated and rotated inner products.
    The setting is that all $\nim^2$ pair-wise alignments must be
    found between $\nim$ images and $\nim$ templates.
    For $\nim\gg 1$, precomputations that scale as $\nim$ are negligible, so
    we give in the second column only the leading term scaling with $\nim^2$.
    Each image (and template) is $n\times n$ pixels.
    We assume that the desired number of rotation angles is $\bigO(n)$, while the number of translations is $\ntrans$.
    The parameter $\nnodes$ is the total number of SVD terms for all the Bessel orders in the expansion (see Section~\ref{sec:svd}); by Theorem~\ref{t:rank} we have $\nnodes = \bigO((W+\log(1/\epsilon))^2)$, where $2W$ is the maximum translation magnitude in pixels, and $\epsilon$ is the desired accuracy.
\label{table:complexity}
}
\end{table}

\section{Problem setup and Fourier representation}
\label{sec:setup}

We use the vectors $\vx, \vd \in \Real^2$ to represent positions and displacements in the spatial domain, while $\vk \in \Real^2$ represents spatial frequencies.
We define the polar coordinates of these vectors as follows:
\begin{eqnarray}
\vx &=& (x \cos \theta, x \sin \theta) \\
\vd &=& (d \cos \omega, d \sin \omega) \\
\vk &=& (k \cos \psi, k \sin \psi) \text{.}
\end{eqnarray}
The Fourier transform of a two-dimensional function $A \in \Leb^2(\Real^2)$ is defined as
\begin{equation}
  \fourier{A}(\vk) := \iint_{\Real^2} A(\vx) \euler^{-\imunit \vk \cdot \vx} \mathop{d\vx}~.
  \label{Ahat}
\end{equation}
The inner product between two functions $A, B \in \Leb^2(\Real^2)$ is written
\begin{equation}
  \langle A, B\rangle = \iint_{\Real^2} A(\vx) B(\vx)^\ast \mathop{d\vx}
  ~,
\end{equation}
where $z^\ast$ is the complex conjugate of $z \in \Complex$.
We will also need Plancherel's theorem \cite{bracewell},
\begin{equation}
  \langle A, B\rangle = \frac{1}{(2\pi)^2} \langle\fourier{A},\fourier{B}\rangle
  ~, \qquad \forall A, B \in \Leb^2(\Real^2)~.
  \label{planch}
\end{equation}

We represent any image as a function $A \in \Leb^2(\Real^2)$, with values corresponding to the image intensity at each location.
We recover $A$ from $\fourier{A}$ using the inverse Fourier transform \cite{bracewell}:
\begin{equation}
A(\vx) = \frac{1}{(2\pi)^2} \iint_{\Real^2} \fourier{A}(\vk) \euler^{\imunit \vk \cdot \vx} \mathop{d\vk} \text{.}
\end{equation}

Let us denote the disk of radius $R>0$ by
\begin{equation}
\Omega_R := \{ \vx \in \Real^2 \mbox{~such~that~} \|\vx\| \leq R \}~.
\end{equation}
We shall assume that all images are supported in $\Omega_1$.
This ensures that $\fourier{A}$ is smooth on the $\bigO(1)$ frequency
scale,
enabling its approximate reconstruction from its values on grids with spacing $\bigO(1)$.
We furthermore assume that $\fourier{A}$ is concentrated in a disk of radius $\kmax$, so that the inversion
\begin{equation}
A(\vx) \approx \frac{1}{(2\pi)^2} \iint_{\Omega_\kmax} \fourier{A}(\vk)\euler^{\imunit \vk \cdot \vx} \mathop{d\vk}
\end{equation}
holds to desired working accuracy.
These assumptions will ensure that we do not lose too much accuracy transitioning to a discrete setting.
\begin{remark}   
  The maximum frequency, $\kmax$, is related to the number of pixels across the image, $n$, as follows.
  The effective image support $\Omega_1$ lives in $[-1,1]$, so that the
  pixel spacing is $\dx = 2/n$, giving the Nyquist
  spatial frequency $\kmax = \pi/\dx = (\pi/2) n$.
  For cryo-EM applications, $n$ is typically $100$ to $400$.
  \label{r:nyq}
\end{remark}   
\begin{remark}   
  In this paper, for simplicity, will assume that $\kmax$ is always the Nyquist frequency,
  so that scalings in $\kmax$ and $n$ are equivalent.
  There do exist applications where $\kmax$ may be much lower
  than $(\pi/2)n$, such as frequency marching \cite{Barnett2017}.
  \label{r:freqmarch}
\end{remark}   
  
  We shall abuse notation slightly and write $A$ and $\fourier{A}$ in polar coordinates as
\begin{eqnarray}
A(x, \theta) &:=& A(\vx) = A(x \cos \theta, x \sin \theta) \\
\fourier{A}(k, \psi) &:=& \fourier{A}(\vk) = \fourier{A}(k \cos \psi, k \sin \psi) \text{.}
\end{eqnarray}
With these reparametrizations, we may now define the rotation of $A$ by $\gamma \in [0, 2\pi)$ as
\begin{equation}
\label{eq:rot-def}
R_\gamma A(x, \theta) := A(x, \theta-\gamma) \text{.}
\end{equation}
Since rotations commute with the Fourier transform, we have the same relationship in the frequency domain:
\begin{equation}
\wfourier{R_\gamma A}(k, \psi) = R_\gamma \fourier{A}(k, \psi) = \fourier{A}(k, \psi-\gamma) \text{.}
\end{equation}

Let us similarly define translation of an image $A$ by the shift vector $\vd \in \Real^2$ as
\begin{equation}
T_\vd A(\vx) := A(\vx - \vd) \text{.}
\end{equation}
In the Fourier domain, abusing notation slightly, we define the action of $T_\vd$ on the Fourier transform $\fA$ as giving the Fourier transform of $T_\vd A$:
\begin{equation}
  T_\vd \fourier{A}(\vk) := \wfourier{T_\vd A}(\vk) = F(\vd, \vk) \fourier{A}(\vk) \text{,}
  \label{Tdeltahat}
\end{equation}
where the \emph{translation kernel} $F(\vd,\vk)$ is defined as
\begin{equation}
  F(\vd,\vk) := \euler^{-\imunit \vd \cdot \vk} = \euler^{-\imunit \delta k \cos(\psi-\omega)} \text{,}
  \label{Fdk}
\end{equation}
which is a plane wave in Fourier space.
We introduce $W$ as a dimensionless user-defined parameter
that specifies the maximum translation magnitude
considered, in wavelength units at the largest frequency $\kmax$.
When $\kmax$ is the Nyquist frequency, this wavelength is $2\dx$.
We thus have translations $\vd$ of magnitude $\delta\le\dmax$ (i.e., $\vd \in \Omega_\dmax$), where $\dmax = 2\pi W/\kmax = 2\dx W$.

Let $B \in \Leb^2(\Real^2)$ be the image against which $A$ is to be aligned
with respect to rotations and translations.
The desired inner product
  $\iint_{\Real^2} R_\gamma T_\vd A(\vx) \, B(\vx)^\ast \mathop{d\vx}$
is, up to a prefactor (which we drop),
well approximated by
\begin{equation}
\hspace{-10mm}
  {\cal X}(\vd, \gamma) \;:=\;
  \iint_{\Omega_\kmax} R_\gamma T_\vd \fourier{A}(\vk) \, \fourier{B}(\vk)^\ast \mathop{d\vk} = \iint_{\Omega_\kmax} F(\vd, \vk) \, \fA(\vk) \, R_{-\gamma} \fB(\vk) \mathop{d\vk}~,
  \label{Xdg}
\end{equation}
which follows from \eqref{planch} and our assumption that $\fourier{A}$ and $\fourier{B}$ are concentrated in $\Omega_\kmax$.
We call \eqref{Xdg} the \textit{bandlimited inner product}.

\subsection{Conversion of discrete input images to the Fourier domain}
\label{sec:pixels}

The above formulation is described in terms of continuous images $A$ and $B$.
However, we are typically provided images sampled on a discrete grid of size $n \times n$ within $[-1, 1]^2$.
Let us consider the pixel grid as representing the set of point samples
\begin{equation}
  \hspace{-5ex}
A_{ij} = A(\vx_{ij})~,\qquad \vx_{ij} := (i\dx -1, j\dx-1)~, \quad i, j \in \{0, \ldots, n-1\}
\end{equation}
where $\dx = 2/n$ is the pixel spacing.
We approximate its Fourier transform $\fA$ at any $\vk \in \Real^2$
by the trapezoid quadrature rule
\begin{equation}
  \fA(\vk) = (\dx)^2 \sum_{i,j=1}^n A_{ij} \euler^{-\imunit \vk \cdot \vx_{ij}}~.
  \label{Ahattrap}
\end{equation}
This quadrature is expected to be accurate if
the image is sufficiently well sampled at the resolution $\dx$
(i.e., before sampling, its frequency content above the Nyquist frequency was negligible),
and $\|\vk\|\le \kmax$.
We similarly define $\fB(\vk)$ given samples $B_{ij}$ for $i, j \in \{0, \ldots, n-1\}$.
We will need to evaluate $\fourier{A}(\vk)$ for $\vk$ on various grids,
which will be discussed later.

Plugging in these Fourier transforms $\fA$ and $\fB$ into \eqref{Xdg}, we obtain a set of inner products ${\cal X}(\vd, \gamma)$ describing the alignment of the discrete images.
The goal of this work is to provide an efficient and accurate method of calculating ${\cal X}(\vd, \gamma)$ from given pixel values $A_{ij}$, $B_{ij}$,
for a large number of user-specified values of $\vd \in \Omega_\dmax$ and $\gamma \in [0, 2\pi)$.

\section{Fourier--Bessel decomposition}
\label{sec:bessel}

In this section we derive a formula for the inner product ${\cal X}(\vd, \gamma)$ defined in
\eqref{Xdg} that separates the dependence on $\vd$ from the dependence on $\vk$. For
this we first need to introduce the Fourier--Bessel basis representation
of each image. At the end of the section we describe how these
formulae are evaluated in practice given $n\times n$ pixel images.

For fixed $k$, $\fA(k, \psi)$ is a $2\pi$-periodic function of $\psi$.
We may therefore decompose the function $\psi \mapsto \fA(k, \psi)$ as a Fourier series, obtaining
\begin{equation}
\fA(k, \psi) = \sum_{q=-\infty}^\infty \bA(k;q) \euler^{\imunit q \psi} \text{,}
\end{equation}
Its coefficients $\bA(k;q)$ are given by
\begin{equation}
  \bA(k;q) = \frac{1}{2\pi} \int_0^{2\pi} \fA(k, \psi) \euler^{-\imunit q \psi} \mathop{d\psi}~.
  \label{aqk}
\end{equation}
By substituting \eqref{Ahat}, changing variable to $\alpha = \theta+\pi/2-\psi$,
and using the integral form of the Bessel function \cite[Eq.~(10.9.1)]{dlmf},
\begin{equation}
  J_n(z) = \frac{1}{2\pi}
  \int_0^{2\pi}
  \euler^{\imunit(z\sin \alpha - n\alpha)} d\alpha~,
  \label{Jint}
\end{equation}
it is easily verified that 
the coefficients are, in terms of the original spatial image,
\begin{equation}
  \bA(k;q) =   \int_0^{2\pi} \int_0^1
  A(x, \theta) \, \left( \euler^{\imunit q(\theta+\pi/2)} J_q(kx) \right)^\ast \, \mathop{x dx d\theta}~.
  \label{aqkbes}
\end{equation}
This takes the form of an inner product of the image $A$ with the 2D
``Fourier--Bessel''
function $u(x,\theta) = \euler^{\imunit q(\theta+\pi/2)} J_q(kx)$, being a polar separation of variables
solution to the Helmholtz equation $(\Delta + k^2)u=0$ at frequency $k$,
where $\Delta$ is the 2D Laplace operator.
This justifies our naming \eqref{aqk} the \emph{Fourier--Bessel coefficients}
of the image.
A similar basis is used by Zhao {\it et al.}\ \cite{Zhao2014,Zhao2016}, who make different assumptions regarding the support of the images.

\begin{remark}[Decay of Bessel functions]   
The large-order asymptotics of $J_q(z)$ may be summarized as follows.
For any fixed $z$, the formula
$J_q(z) \sim (2\pi q)^{-1/2} (\euler z/2 q)^q$ shows super-exponential
decay in $q$ \cite[Eq.~(10.19.1)]{dlmf}.
More precisely, for fixed real $z$, $J_q(z)$ is
oscillatory and significant in the order range $q \le z$,
and exponentially small for $q \ge z + \bigO(z^{1/3})$,
the latter statement following from the Debye
asymptotics \cite[Eq.~(10.19.8)]{dlmf}.
Loosely speaking, this is associated
with the ``bump'' of the Bessel function
(as a function of $z$) centered at
$z\approx q$ having width roughly $z^{1/3} \approx q^{1/3}$.
\label{r:Jdecay}
  \end{remark}     

The above remark, \eqref{aqkbes}, and $x\le1$ show that the coefficients
$\bA(k;q)$ are exponentially small, and thus numerically negligible,
for $|q| > K + \bigO(K^{1/3})$.

We now derive the actions of rotation and translation on these
Fourier--Bessel coefficients.
Applying \eqref{aqk} to the rotated Fourier transform $R_\gamma \fA$, we obtain
\begin{eqnarray}
\nonumber
&& \frac{1}{2\pi} \int_0^{2\pi} R_\gamma \fA(k, \psi) \euler^{-\imunit q \psi} \, d\psi
= \frac{1}{2\pi} \int_0^{2\pi} \fA(k, \psi-\gamma) \euler^{-\imunit q \psi} \, d\psi  \\
&=& \frac{1}{2\pi} \int_0^{2\pi} \fA(k, \psi) \euler^{-\imunit q (\psi + \gamma)} \, d\psi = \euler^{-\imunit q \gamma} \bA(k; q) \text{.}
\end{eqnarray}
This shows that a rotation acts on the Fourier--Bessel coefficients $\bA$
as a diagonal multiplication operator,
\begin{equation}
R_\gamma \bA(k, q) := \euler^{-\imunit q \gamma} \, \bA(k; q)~.
\end{equation}
Here we abuse notation slightly by using the same symbol
$R_\gamma$ that we used for rotations in an angular variable.
This diagonal action will allow us to simultaneously compute ${\cal X}(\vd, \gamma)$ for multiple values of $\gamma$ in an efficient manner.

Let us similarly denote by $T_\vd a(k;q)$ the action of translation by the vector
$\vd$ on the Fourier--Bessel coefficients $\bA$.
Using \eqref{Tdeltahat}, \eqref{Fdk} and \eqref{aqk} we compute
\begin{eqnarray}
\hspace{-10ex}T_\vd a(k;q)
  &=& \frac{1}{2\pi} \int_0^{2\pi} T_\vd \fA(k, \psi) \euler^{-\imunit q \psi} \, d\psi
= \frac{1}{2\pi} \int_0^{2\pi} \fA(k, \psi) \euler^{-\imunit k \delta \cos(\psi-\omega)} \euler^{-\imunit q \psi} \, d\psi \nonumber \\
&=& \frac{1}{2\pi} \int_0^{2\pi} \fA(k, \psi) \left( \sum_{\ell=-\infty}^\infty J_\ell(\delta k) \euler^{\imunit \ell (\psi - \omega - \pi/2)} \right) \euler^{-\imunit q \psi} \, d\psi \nonumber \\
&=& \sum_{\ell=-\infty}^\infty J_\ell(\delta k) \euler^{-\imunit \ell (\omega + \pi/2)} \, \bA(k; q-\ell) \text{,}
\end{eqnarray}
where on the third line we
used the Jacobi--Anger expansion \cite[10.12.3]{dlmf},
\begin{equation}
\label{eq:jacobi-anger}
\euler^{-\imunit z \cos \phi} = \sum_{\ell=-\infty}^{\infty} J_\ell(z) \euler^{\imunit \ell (\phi - \pi/2)}~.
\end{equation}
In summary, translation acts as a 1D discrete convolution on the
coefficients,
\begin{equation}
\label{eq:fbk1}
T_\vd \bA(k, q) := \sum_{\ell=-\infty}^\infty f(\vd, k; \ell) \bA(k; q-\ell) \text{,}
\end{equation}
where $f(\vd, k; \ell)$ denotes the Fourier--Bessel translation kernel
\begin{equation}
  f(\vd, k; \ell) \;:=\;
  J_\ell(\delta k) \euler^{-\imunit \ell (\omega + \pi/2)}~,
  \label{fdef}
\end{equation}
being the $\ell$th 1D Fourier series coefficient of $F(\vd,\vk)$
on the ring $\|\vk\|=k$ (see \eqref{Fdk}).
Note that the number of significant terms in the convolution is not large,
and scales only like $W$, the size of the maximum translation in wavelengths.
This follows from Remark~\ref{r:Jdecay}, and \eqref{fdef},
where the maximum Bessel argument of $\dmax\kmax = 2\pi W$
shows that $f_\ell(\vd,k)$
is exponentially small for $|\ell| > 2\pi W + \bigO((2\pi W)^{1/3})$.

It will turn out that each of the terms in the above convolution is amenable to factorization, which will help us consider multiple translations $\vect{\delta}$.
Indeed the factors dependent on $\vect{\delta}$ ($\delta$ and $\omega$) and those dependent on $\vect{k}$ ($k$ and $\psi$) are \emph{almost} separable from one another.
The exception is the Bessel term $J_{\ell}(k \delta)$, which we return to in Section~\ref{sec:svd}.

Now we turn to the desired inner product.
From \eqref{aqk}, on each fixed $k$ ring we have the 1D Plancherel formula
\begin{equation}
\int_0^{2\pi} \fourier{A}(k,\psi)\fourier{B}(k,\psi)^\ast d\psi = 2\pi \sum_{q=-\infty}^\infty a(k;q)b(k;q)^\ast~.
\end{equation}
Using this we rewrite the inner product \eqref{Xdg} as
\begin{eqnarray}
\hspace{-10mm}
\nonumber
{\cal X}(\vd, \gamma) &=& \iint_{\Omega_\kmax} R_\gamma T_\vd \fourier{A}(\vk) \, \fourier{B}(\vk)^\ast \mathop{d\vk}
= 2\pi \sum_{q=-\infty}^\infty \int_0^\kmax R_\gamma T_\vd \bA(k; q) \, \bB(k; q)^\ast \mathop{kdk} \\
\nonumber
&=& 2\pi \sum_{q=-\infty}^\infty \int_0^\kmax \euler^{-\imunit q \gamma} \left( \sum_{\ell=-\infty}^\infty f(\vd, k; \ell) \bA (k; q-\ell) \right) \bB(k; q)^\ast \mathop{kdk} \\
\label{eq:bessel-inner}
&=& 2\pi \sum_{q=-\infty}^\infty \sum_{\ell=-\infty}^\infty \euler^{-\imunit q \gamma} \int_0^\kmax f(\vd, k; \ell) \bA(k; q-\ell) \bB(k; q)^\ast \mathop{kdk}
\end{eqnarray}
While this expression appears more complicated than our original formula \eqref{Xdg}, we shall see that it lends itself more naturally to a fast numerical algorithm.
For now we will simply remark that this expression displays the \emph{trilinearity} of the inner product: ${\cal X}(\vd, \gamma)$ is linear in $\fA(\vk)$, $\fB(\vk)$ and $f(\vd, k; \ell)$.
Our method will exploit this by reducing the approximation of ${\cal X}(\vd, \gamma)$ to an approximation of the translation kernel $F(\vd, k; \ell)$.

\subsection{Discretization and truncation}
\label{sec:discretization}

To numerically approximate ${\cal X}(\vd, \gamma)$ evaluated via
\eqref{eq:bessel-inner}, we make the following approximations.
First, we truncate the infinite sums over $q$ and $\ell$ to finite sums over $|q| \leq \qmax$ and $|\ell| \leq \ellmax$.
Second, we replace the integral over $k$ with a quadrature scheme designed on $[0, K]$.
We prefer Gauss--Jacobi quadrature built with a weight function corresponding to the area-element $kdk$.
This quadrature scheme uses Jacobi polynomials of type $\alpha=0$, $\beta=1$ (i.e.\ corresponding to weight $k$ after a change of variable), and produces a set of nodes $k_m \in [0, K]$ and weights $w_m$ for $m \in \{1, \ldots, \mmax\}$ such that
\begin{equation}
  \int_{0}^{\kmax} g(k) \mathop{kdk}
  \;\approx\; \sum_{m=1}^{\mmax} g(k_m)w_m
\end{equation}
holds to high accuracy for any radial function $g$ that is smooth
on the scale of $\bigO(1)$.
In summary, we approximate ${\cal X}(\vd, \gamma)$ by
\begin{equation}
\hspace{-15mm}
X(\vd, \gamma) = 2\pi \sum_{q=-\qmax}^\qmax \sum_{\ell=-\ellmax}^\ellmax \sum_{m=1}^\mmax \euler^{-\imunit q \gamma} \bA(k_m; q-\ell) f(\vd, k_m; \ell) \bB(k_m; q)^\ast w_m \text{.}
\label{Xtrunc}
\end{equation}
We will not attempt a rigorous convergence analysis here,
but rather justify the empirical sizes of the convergence parameters
$\qmax$, $\ellmax$, and $\mmax$ needed.
As discussed in the previous section, following Remark~\ref{r:Jdecay},
$\qmax$ need only slightly exceed $\kmax=(\pi/2)n$.
Similarly, $\ellmax$ need only slightly exceed $2\pi W$.
For quadrature, in order to capture the $\bigO(\kmax)$ oscillations in
the integrand, we need $\mmax = \bigO(\kmax) = \bigO(n)$.
Once $M$ is large enough to resolve these oscillations, the
convergence becomes super-algebraic;
this follows because
the integrands (deriving from Fourier transforms of images with compact support)
are analytic.
It is thus easy to choose the parameters to achieve a desired accuracy (see Section~\ref{sec:results}).

\begin{remark}
We could reduce the number of operations in the above calculation by making $\qmax$ depend on $k$.
In particular, we expect $\qmax$ to vary proportionally to $k$, as it represents the number of angular modes a ring of circumference $2\pi k$.
This would reduce the total number of Fourier--Bessel coefficients by a factor of two.
However, this reduction does not necessarily translate to shorter running times, as the resulting triangular sum cannot easily take advantage of vectorized libraries for FFTs and matrix multiplication.
\end{remark}

\subsection{Decomposition of discrete images}
\label{sec:nufft}

The approximate inner product \eqref{Xtrunc}
is defined using the Fourier--Bessel coefficients $\bA(k_m,q)$ and $\bB(k_m,q)$
on all radial quadrature nodes $\{k_m\}_{m=1}^\mmax$, for all indices
$|q|\le \qmax$. This is $\bigO(n^2)$ coefficients.
Given $n\times n$ pixel images $\{A_{ij}\}$ and $\{B_{ij}\}$
as in Section~\ref{sec:pixels},
we compute these coefficients
in optimal complexity as follows.
For each image we apply the
2D type 2 non-uniform FFT (NUFFT) \cite{dutt,nufft,usingnfft}
(we use the FINUFFT implementation \cite{finufft})
to evaluate the Fourier transform \eqref{Ahattrap} on
a tensor-product polar $\vk$
grid $\{(k_m,\psi_p)\}_{m=1,\dots,\mmax, \; p=1,\dots,2\qmax}$
where the uniform angular grid has $\psi_p := \pi (p-1)/\qmax$.
The NUFFT cost depends weakly on a tolerance parameter $\epsilon$, which we
ensure is smaller than the desired overall accuracy.
The complexity of this step is $\bigO(n^2 \log n + n^2 \log (1/\epsilon))$.

Finally, on each quadrature ring $m=1,\dots,\mmax$,
we approximate \eqref{aqk} using the $\fourier{A}$ values on that ring
with the periodic trapezoid rule,
\begin{equation}
  \bA(k;q) = \frac{1}{2\pi} \int_0^{2\pi} \fA(k, \psi) \euler^{-\imunit q \psi} \mathop{d\psi}
  \approx \frac{1}{2\qmax} \sum_{p=1}^{2\qmax} \fA(k, \psi_p) \euler^{-\imunit q \psi_p}~.
  \label{aqkcalc}
\end{equation}
  Since $\psi_p = \pi (p-1)/\qmax$, this takes the form of a 1D discrete Fourier
  transform of length $2\qmax$, which we evaluate using the FFT.
The total cost of these $\bigO(n)$ 1D FFTs is $\bigO(n^2 \log n)$.
In summary, fixing $\epsilon$, the
entire precomputation then scales as $\bigO(n^2 \log n)$ per image.

\section{Fast computation of multiple inner products}
\label{sec:fast}

In this section, we present an SVD interpolation
method for computing many translated
and rotated inner products \eqref{Xtrunc}.
For simplicity we consider one image $A$ and one template $B$.
To provide intuition,
we start with an algorithm in which translations are treated
independently,
explain interpolation with respect to translations,
then explain our method as an optimal version of such an interpolation.

\subsection{Multiple rotation angles}
\label{sec:fastrots}

We first present a variant of the BFT algorithm.
If we only have a single shift vector $\vd$, we can calculate $X(\vd, \gamma)$ efficiently for a large number of $\gamma$ values by taking advantage of the diagonal action of $R_\gamma$ in the Fourier--Bessel coefficients.
Indeed, we have
\begin{equation}
\label{eq:fastrots}
\hspace{-15mm}
X(\vd, \gamma) = 2\pi \sum_{q=-\qmax}^\qmax \euler^{-\imunit q \gamma} \sum_{m=1}^\mmax w_m \bB(k_m; q)^* \sum_{\ell=-\ellmax}^\ellmax \bA(k_m; q-\ell) f(\vd, k_m; \ell) \text{.}
\end{equation}
This calculation can be achieved in two steps:
\begin{eqnarray}
\label{eq:fbk2}
&& \textbf{[Step 1]} \quad \fourier{X}(\vd, q) = 2\pi \sum_{m=1}^\mmax w_m \bB(k_m; q)^* \sum_{\ell=-\ellmax}^\ellmax \bA(k_m; q-\ell) f(\vd, k_m; \ell) \\
&& \textbf{[Step 2]} \quad X(\vd, \gamma) = \sum_{q=-\qmax}^\qmax \euler^{-\imunit q \gamma} \fourier{X}(\vd, q) \label{1dfft} \text{.}
\end{eqnarray}
The first step computes the Fourier coefficients $\fourier{X}(\vd, q)$, requiring $\bigO(\mmax \qmax \ellmax)$ operations.
If we restrict $\gamma$ to the set $\{0, \pi/\qmax, \ldots, 2\pi (2\qmax-1) / \qmax\}$, we may evaluate the second step using a 1D FFT of size $2\qmax$.
The first step dominates,
so the total complexity is $\bigO(n^2 \wmax)$, where we have used $\qmax = \bigO(n)$, $\mmax = \bigO(n)$ and $\ellmax = \bigO(\wmax)$ (see previous section).
By zero-padding and using a larger FFT (or NUFFT) in the second step,
a large number of $\gamma$ values may be sampled essentially for free.
A naive method to now handle $\ntrans$ different translations
$\vd \in \Omega_\dmax$ is to repeat the above two steps, separately for each
translation,
giving a cost $\bigO(\ntrans n^2 \wmax)$.
For each image-template pair, this is $W$ times slower than BFT;
see Table~\ref{table:complexity}.
However, the above steps,
with $f(\vd, k_m; \ell)$ replaced by other functions, will serve as a building block for the much faster algorithm described in the following subsections.

\subsection{Linear interpolation over translations} 
\label{sec:linear}

We now present a simple way to accelerate the brute force inner product calculation over multiple translations.
Instead of evaluating $X(\vd, \gamma)$ at all values of $\vd$ which are of interest to us, one can use the fact that, for fixed $\gamma$, the mapping $\vd \mapsto X(\vd, \gamma)$ is smooth.
Indeed, it is the smoothness of the underlying image functions
$A$ and $B$ induces smoothness in $X$.
A common approach is to calculate $X(\vd, \gamma)$ at some subset of shifts and use local (low-order)
interpolation to obtain its value at other locations.

Let us denote such a set of $\nnodes$ interpolation nodes by $\{ \vd_1, \ldots, \vd_{\nnodes} \} \subset \Omega_\dmax$.
This could be a grid (e.g., Cartesian or polar) covering $\Omega_\dmax$;
we assume it is independent of $\gamma$.
Fixing $\gamma$ for now,
we first compute $X(\vd_\zeta, \gamma)$ for $\zeta \in \{1, \ldots, \nnodes\}$,
using, for instance the method outlined in the previous section.
For other values of $\vd$ at this $\gamma$, we interpolate
\begin{equation}
\label{lininterp}
X^{\lin}(\vd, \gamma) = \sum_{\zeta=1}^{\nnodes} Y_\zeta^\lin(\vd) X(\vd_\zeta, \gamma) \text{,}
\end{equation}
where $Y_\zeta^\lin(\vd)$ is the interpolation weight (basis function) associated with $\vd$ and $\vd_\zeta$.
In the bilinear case, it is nonzero only for the four nodes $\vd_\zeta$
that are nearest to $\vd$.
The above is repeated for the $\bigO(\kmax)$ desired $\gamma$ values.
Since the complexity for computing $X(\vd_\zeta, \gamma)$ is $\bigO(n^2 \wmax)$ (see Section~\ref{sec:fastrots}), the overall complexity is $\bigO(\nnodes n^2 \wmax)$ for all $\zeta$ (this drops to $\bigO(\nnodes n^2)$ if we instead use the BFT, see Table~\ref{table:complexity} and Section~\ref{sec:bft}).
The cost of local interpolation is $\bigO(1)$ per
lookup, or $\bigO(\ntrans n)$ per image-template pair.
This is advantageous since the original $\ntrans n^2$ factor is gone.
However, to achieve high interpolation accuracy in a low-order scheme,
a fine grid of nodes would be needed, hence a large $\nnodes$.

\subsection{Least-squares interpolation}
\label{sec:lsq}

We may now ask: can we increase the accuracy without significantly increasing the computational cost?
First note that the interpolated result can be rewritten by moving the sum over $\zeta$ inside the other sums:
\begin{equation}
\hspace{-20mm}
X^{\lin}(\vd, \gamma) = 2\pi \sum_{q=-\qmax}^\qmax \euler^{-\imunit q \gamma} \sum_{m=1}^\mmax w_m \bB(k_m; q)^* \sum_{\ell=-\ellmax}^\ellmax \bA(k_m; q-\ell) \sum_{\zeta=1}^{\nnodes} Y_\zeta(\vd) f(\vd_\zeta, k_m; \ell) \text{.}
\label{Xlsq}
\end{equation}
Comparing this to the formula for $X(\vd, \gamma)$ in \eqref{eq:fastrots}, we see that linear interpolation is equivalent to the approximation
\begin{equation}
  \sum_{\zeta=1}^{\nnodes} Y_\zeta^\lin(\vd) f(\vd_\zeta, k; \ell)
  \;\approx \;
  f(\vd, k; \ell) \text{,} \quad \forall k \in [0, \kmax], \ell \in \Integer \text{.}
\end{equation}
In other words, due to the linearity of the bandlimited inner product $X(\vd, \gamma)$, an interpolation scheme over $X(\vd, \gamma)$ corresponds to an interpolation scheme over $f(\vd, k; \ell)$.

Fixing $\vd$, we now present a scheme to
find the vector of interpolation weights
$Y^\ls(\vd) := (Y_1^\ls(\vd), \ldots, Y_{\nnodes}^\ls(\vd))$
that best approximates $f(\vd, k; \ell)$ in a least-squares sense.
Equation \eqref{Xlsq} suggests that one should
minimize the sum-of-squares error over the quadrature nodes $k_1, \ldots, k_\mmax$ and angular modes $-\ellmax \le \ell \le \ellmax$;
however, treating $k$ as a continuous variable is cleaner, gives
essentially the same answer (since the quadrature scheme is accurate),
and will allow analysis relevant for the SVD method.

The weights vector is the least-squares solution
\begin{equation}
\hspace{-15mm}
Y^\ls(\vd) = \argmin_{Y=(Y_1, \ldots, Y_{\nnodes})} \sum_{\ell=-\infty}^\infty \int_0^\kmax \left| \sum_{\zeta=1}^{\nnodes} Y_\zeta f(\vd_\zeta, k; \ell) - f(\vd, k; \ell) \right|^2 \, kdk~.
\end{equation}
The corresponding normal equations for $Y^\ls$ form the $\nnodes \times \nnodes$ linear system
\begin{equation}
\label{eq:normal2d}
\sum_{\zeta=1}^{\nnodes} \besselker^\ls(\vd_{\zeta^\prime}, \vd_\zeta) Y_\zeta^\ls(\vd) = \besselker^\ls(\vd_{\zeta^\prime}, \vd)~,\qquad \forall \zeta^\prime \in \{1, \ldots, \nnodes\} \text{,}
\end{equation}
where
\begin{equation}
\besselker^\ls(\vd^\prime, \vd) := \sum_{\ell=-\infty}^\infty \euler^{\imunit \ell (\omega^\prime - \omega)} \int_0^\kmax J_\ell(\delta^\prime k) J_\ell(\delta k) \mathop{kdk} \text{.}
\end{equation}
These integrals may be calculated analytically:
\begin{equation}
\label{eq:normal2dA}
\besselker^\ls(\vd^\prime, \vd) = J_{0}(\tilde{\delta}\kmax) + J_{1}(\tilde{\delta}\kmax)~\text{,}
\end{equation}
where $\tilde{\delta} = |\vd^\prime-\vd|$.
The resulting linear system is well-conditioned for moderate values of $\dmax$ and $\kmax$ (and hence $\wmax$), so $Y_\zeta^\ls(\vd)$ may be calculated to high precision using standard numerical linear algebra techniques.
This is an image-independent precomputation required for each $\vd$.
Given the same set of $X(\vd_{\zeta},\gamma)$ samples \eqref{lininterp}, the interpolator is now
\begin{equation}
X^{\ls}(\vd, \gamma) = \sum_{\zeta=1}^{\nnodes} Y_\zeta^\ls(\vd) X(\vd_\zeta, \gamma)~.
\end{equation}
This increases the lookup cost from $\bigO(1)$ per shift to $\bigO(\nnodes)$,
but achieves higher accuracy.

We can improve the accuracy by optimizing the node locations $\{\vd_1, \ldots, \vd_{\nnodessls}\}$ using, for example, gradient descent.
The gradient is readily calculated by noting that
\begin{equation}
\label{eq:normal2dB}
\partial_{\tilde{\delta}\kmax} \besselker^\ls(\vd^\prime, \vd) = \frac{1}{2}J_{-1}(\tilde{\delta}\kmax) - \frac{1}{2}J_{3}(\tilde{\delta}\kmax)~\text{.}
\end{equation}
While this works well for a small number of nodes, this nonconvex optimization problem becomes intractable for larger $\nnodessls$,
limiting applicability to lower accuracies ($\epsilon \ge 10^{-2}$).

\subsection{Generalized least-squares interpolation}
\label{sec:slsq}

One could potentially increase the accuracy of the above
method by allowing a {\em different} vector of interpolation weights
$Y_\zeta^\sls(\vd; \ell) := (Y_{\zeta}^{\sls}(\vd;\ell))_{\zeta=1}^\nnodessls$
for each mode $\ell \in \{-\ellmax,\dots,\ellmax\}$.
We refer to this as the \textit{generalized least-squares interpolation} (GLS) method.
We now solve
\begin{equation}
\hspace{-10mm}
Y^\sls(\vd; \ell) = \argmin_{Y=(Y_1, \ldots, Y_{\nnodessls})} \int_0^\kmax \left| \sum_{\zeta=1}^{\nnodessls} Y_\zeta f(\vd_\zeta, k; \ell) - f(\vd, k; \ell) \right|^2 \mathop{kdk} \text{,}
\end{equation}
Using \eqref{fdef} we see that $Y_{\zeta}^\sls(\vd; \ell) = \euler^{-\imunit \ell (\omega - \omega_\zeta)} U^\sls_\zeta(\delta; \ell)$, where
\begin{equation}
U^\sls_\zeta(\delta; \ell) = \argmin_{U=(U_1, \ldots, U_\nnodessls)} \int_0^\kmax \left| \sum_{\zeta=1}^\nnodessls U_\zeta J_\ell(\delta_\zeta k) - J_\ell(\delta k) \right|^2 \mathop{kdk} \text{.}
\end{equation}
As before, $U^\sls_\zeta(\delta; \ell)$ satisfies the normal equations
\begin{equation}
\label{normalsls}
\sum_{\zeta=1}^{\nnodessls} \besselker^\sls(\delta_{\zeta^\prime}, \delta_\zeta; \ell) U^\sls_\zeta(\delta; \ell) = \besselker^\sls(\delta_{\zeta^\prime}, \delta; \ell) \quad \forall \zeta^\prime \in \{1, \ldots, \nnodessls\} \text{,}
\end{equation}
where
\begin{equation}
\label{besselkersls}
\besselker^\sls(\delta^\prime, \delta; \ell) = \int_0^\kmax J_\ell(\delta^\prime k) J_\ell(\delta k) \mathop{k dk} \text{.}
\end{equation}
As before, \eqref{normalsls} can be readily solved by analytically evaluating \eqref{besselkersls} and applying standard numerical techniques.

Given a set of weights $Y_\zeta^\sls(\vd; \ell)$, we thus have the steps:
\begin{eqnarray}
&& \textbf{[Step 1]} \quad \fourier{Z}^\sls_\zeta(q, \ell) = 2\pi \sum_{m=1}^\mmax w_m \bB(k_m; q)^* \bA(k_m; q-\ell) f(\vd_\zeta, k_m; \ell) \\
&& \textbf{[Step 2]} \quad Z^\sls_\zeta(\gamma; \ell) = \sum_{q=-\qmax}^\qmax \euler^{-\imunit q \gamma} \fourier{Z}^\sls_\zeta(q, \ell) \\
&& \textbf{[Step 3]} \quad X^{\sls}(\vd, \gamma) =  \sum_{\ell=-\ellmax}^\ellmax \sum_{\zeta=1}^{\nnodessls} Y_\zeta^\sls(\vd; \ell) Z_\zeta(\gamma; \ell) \text{.}
\end{eqnarray}
Another advantage over the plain least-squares interpolation method of Section~\ref{sec:lsq}
is that a different $\nnodessls$ could be used for each $\ell$;
for instance, from \eqref{fdef} and Remark~\ref{r:Jdecay} one could
use fewer nodes for large $|\ell|$ where $J_\ell(k\delta)$ is small.
In the next subsection, rather than pursue least-squares methods,
we use this intuition to understand an improved SVD-based method.

\subsection{Proposed SVD interpolation method: \meth}
\label{sec:svd}

In the previous section, we optimized weights $Y_\zeta(\vd; \ell)$ to minimize the mean squared error between $f(\vd, k; \ell)$ and $\sum_{\zeta=1}^\nnodessls Y_\zeta(\vd; \ell) f(\vd_\zeta, k; \ell)$.
To further reduce the error, we need to optimize the second factor, $f(\vd_\zeta, k; \ell)$.
As discussed in Section~\ref{sec:svd}, we could optimize the interpolation nodes, but this is only tractable for low accuracies.

\begin{figure}
\centering
\includegraphics[width=0.70\linewidth,trim= {0cm 0cm 0cm 0cm},clip]{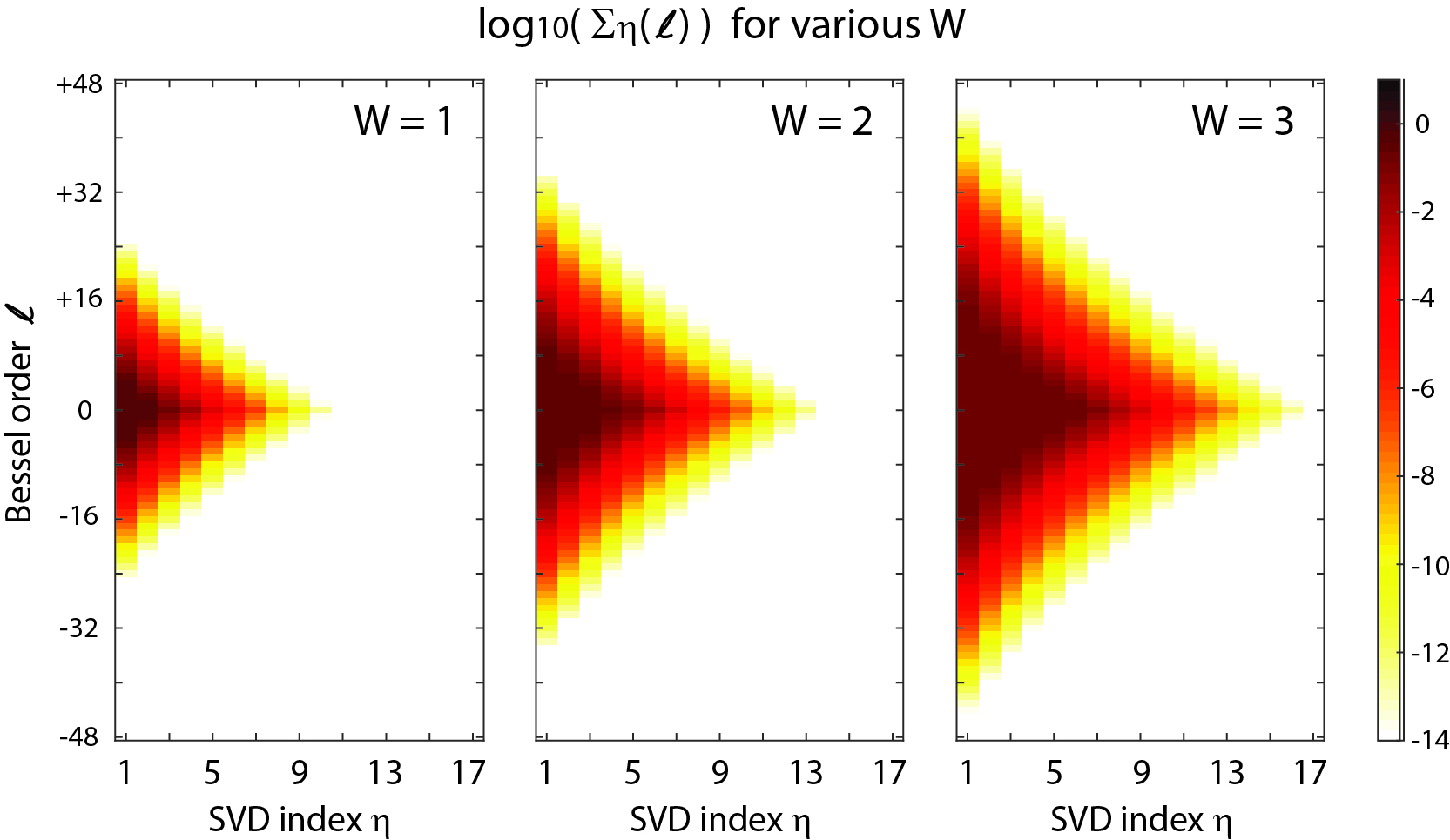}
\caption{\label{fig:Sigma}
Operator singular values $\Sigma_{\eta}(\ell)$ of the Bessel kernel
$J_\ell(k\delta)$ over $\delta\in [0,\dmax]$ and $k\in [0,\kmax]$,
as a function of $\ell$ and singular value index (see Section~\ref{sec:svd}).
The cases $W=1,2$ and $3$ are shown;
$W$ is the maximum translation magnitude in wavelengths ($2\dx$).
In each case the values of $\log_{10}(\Sigma_{\eta}(\ell))$ are displayed using the colorbar to the far right.
The number of SVD terms $\nnodes$ equals the number of index pairs $(\eta,\ell)$
with $\Sigma_\eta(\ell) \ge \epsilon$; these pairs form an approximately
triangular region.
}
\end{figure}

\begin{figure}
\centering
  \includegraphics[width=0.70\linewidth,trim= {0cm 0cm 0cm 0cm},clip]{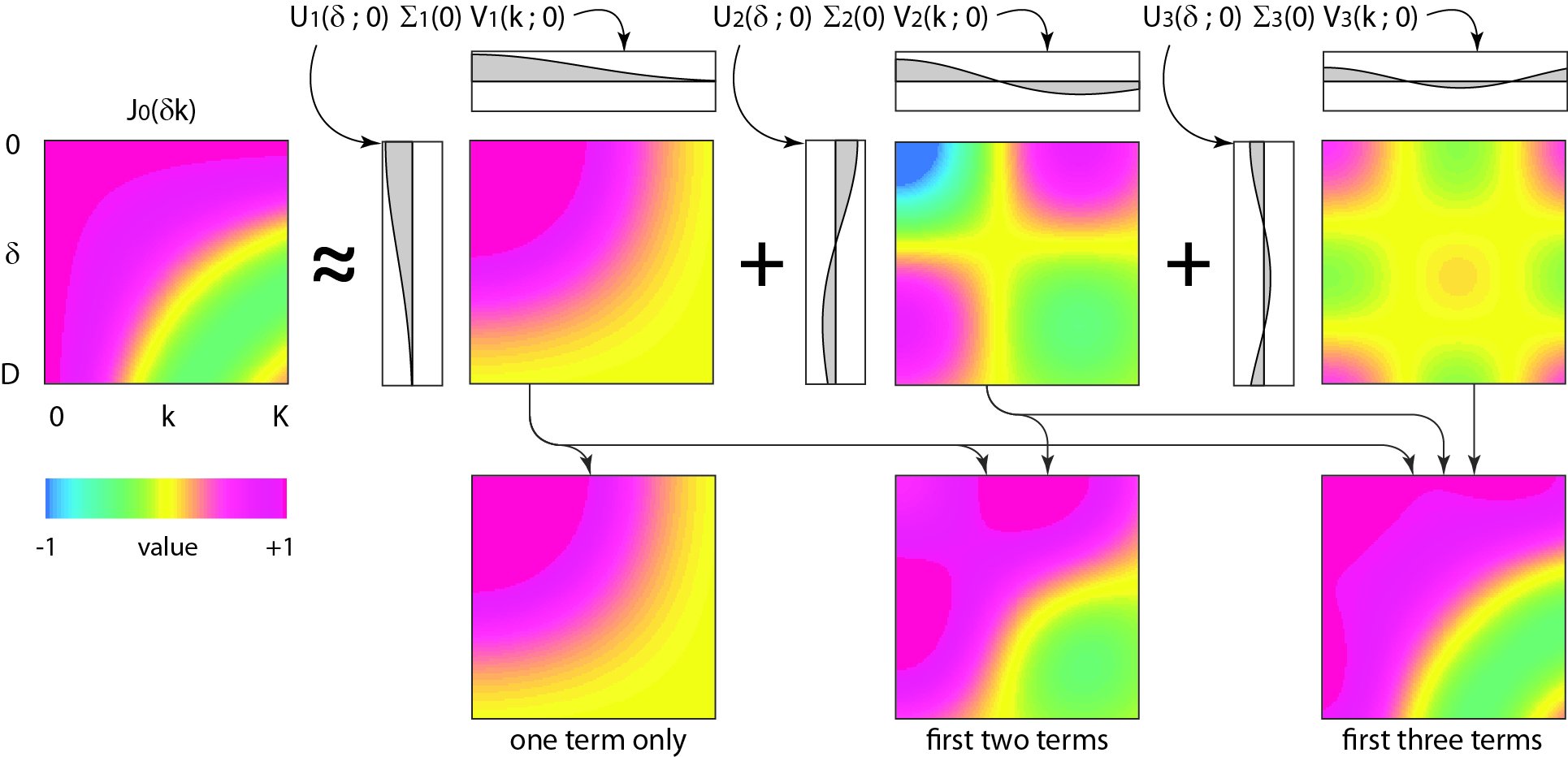}
\caption{\label{fig:bessel-svd}
The image in the top left corner shows the Bessel function $J_0$ as a function of $\delta$ (vertical axis) and $k$ (horizontal axis), i.e.\ the case $\ell=0$.
The domains $\dmax$ and $\kmax$ correspond to $W=1$.
On the left-hand side we show the first three terms in the corresponding SVD.
Each image shows the product $U_\eta(\delta; 0) \Sigma_\eta(0) V_{\eta}(k; 0)$ for $\eta = 1, 2, 3$.
To the left of each image is shown $U_\eta(\delta; 0) \sqrt{\Sigma_\eta(0)}$, while on top we show $\sqrt{\Sigma_\eta(0)} V_\eta(\delta; 0)$.
Below each SVD outer product term, we plot their cumulative sum.
Here $H_{0}=3$ terms is sufficient to capture the coarse features of $J_{0}(\delta k)$.
}
\end{figure}

Treating the interpolant samples
$f(\vd_\zeta, k; \ell)$ as functions of $k$, we now replace them
by a general set of functions $G_\eta(k; \ell)$, where $\eta=1,\dots, \nnodesell$
indexes the $\ell$th set.
Independently for each mode $\ell$,
we choose the set $\{G_\eta(k; \ell)\}_{\eta=1}^{\nnodesell}$ that optimally interpolates
$f(\cdot, k;\ell)$ over the disk $\Omega_\dmax$ in the mean-square sense.
Not only is this optimization problem more general than the least-squares problems above, but, as it turns out, it is also easier to solve.
Thus, for each $\ell$, we minimize the mean squared error, which is the squared \emph{Hilbert--Schmidt} (HS) norm of the residual,
\begin{eqnarray}
\nonumber
  E_\hs(\ell)^2
  &:=& \int_{\Omega_{\dmax}} \int_0^{\kmax} \left| \sum_{\eta=1}^{\nnodesell} Y_\eta(\vd; \ell) G_\eta(k; \ell) - f(\vd, k; \ell) \right|^2 \mathop{kdk\, d\vd}
\\
\nonumber
&=& \int_{\Omega_{\dmax}} \int_0^{\kmax} \left| \sum_{\eta=1}^{\nnodesell} Y_\eta(\vd; \ell) G_\eta(k; \ell) - J_\ell(\delta k) \euler^{-\imunit \ell (\omega + \pi/2)} \right|^2 \mathop{kdk \, d\vd} \\
\nonumber
&=& \int_{\Omega_{\dmax}} \int_0^{\kmax} \left| \sum_{\eta=1}^{\nnodesell} Y_\eta(\vd; \ell) \euler^{\imunit \ell (\omega + \pi/2)} G_\eta(k; \ell) - J_\ell(\delta k) \right|^2 \mathop{kdk \, d\vd} \\
&\label{EFl2} =& 2\pi \int_0^{\dmax} \int_0^{\kmax} \left| \sum_{\eta=1}^{\nnodesell} U_\eta(\delta; \ell) G_\eta(k; \ell) - J_\ell(\delta k) \right|^2 \mathop{k dk \,\delta d\delta} \text{.}
\end{eqnarray}
Here we substituted \eqref{fdef}, then noted that since $J_\ell(\delta k)$ does not depend on the angle $\omega$, the $Y_\eta(\vd; \ell)$ that minimizes the error must be such that the product $Y_\eta(\vd; \ell) \euler^{\imunit \ell (\omega + \pi/2)}$ also does not depend on $\omega$.
We therefore replace this product by $U_\eta(\delta; \ell)$ (which only depends on the magnitude $\delta$) and set $Y_\eta(\vd; \ell) = U_\eta(\delta; \ell) \euler^{-\imunit \ell(\omega + \pi/2)}$.
As a result, we seek a rank-$\nnodesell$ separable approximation of the bivariate function $J_\ell(\delta k): [0, \dmax] \times [0, \kmax] \rightarrow \mathbb{R}$
which is mean-square optimal under the linear weight function.
Treating this bivariate function as the kernel of a (compact) integral operator,
the solution is to truncate to the first $\nnodesell$
terms in the operator SVD \cite[Chap.~II]{GK69}.
As proven by Schmidt in 1907 \cite[\S 18]{schmidt}
(see \cite{stewartFHS}),
the truncated SVD is the optimal fixed-rank approximation to an operator
in the Hilbert--Schmidt norm.
With our weight function, for each $\ell$,
the operator SVD is
\begin{equation}
  J_\ell(\delta k) \;=\; \sum_{\eta=1}^\infty  U_\eta(\delta; \ell) \Sigma_\eta(\ell) V_\eta(k; \ell)~,  \qquad \delta\in[0,\dmax],\quad k\in[0,\kmax]~,
  \label{Jsvd}
\end{equation}
where $\Sigma_1(\ell) \ge \Sigma_2(\ell) \ge \dots \to 0$ are the
singular values, the left singular functions $\{U_\eta(\delta; \ell)\}_{\eta=1}^\infty$ are
orthonormal on $[0,\dmax]$ with respect to $\delta d \delta$,
and the right singular functions $\{V_\eta(k; \ell)\}_{\eta=1}^\infty$ are
orthonormal on $[0,\kmax]$ with respect to $k dk$.
Although the singular values appear to depend on both $\kmax$ and $\dmax$,
one may check by changing variables that in fact they depend
solely on their product $\kmax \dmax=2\pi W$.
Figure~\ref{fig:Sigma} shows the rapid decay of these singular values,
and Figure~\ref{fig:bessel-svd} the first few singular functions.

\begin{remark}[Approximating the operator SVD]  
  In practice, to compute the SVD, we use a tensor-product quadrature scheme
  to project $J_\ell(\delta k)$ onto a product basis of orthogonal Jacobi polynomials with $\alpha = 0$ and $\beta = 1$, rescaled onto the intervals $\delta \in [0, \dmax]$ and $k\in [0, \kmax]$.
  These Jacobi bases are orthonormal with respect to $\delta d\delta$ and $kdk$,
  respectively.
Given the matrix of coefficients representing $J_\ell$ in this basis, we
compute its SVD using standard dense numerical linear algebra \cite{golubvanloan,trefethenbau}.
The left and right singular functions $U_\eta(\delta;\ell)$ and $V_\eta(\delta;\ell)$ are then given by summing the Jacobi bases with
the corresponding singular vector entries as coefficients.
\label{r:svd}
\end{remark}              

\begin{remark}[Alternative weight functions]  
  In this paper we calculate the SVD of each $J_{\ell}$ with the radial weights $\delta d\delta$ and $kdk$, corresponding to
  $L^2$ approximation over the disks $\Omega_{\dmax}$ and $\Omega_{\kmax}$.
One can easily generalize this to accommodate alternative weight functions.
For example, when processing images of known power-spectral decay, one should tailor the $k$-weight to emphasize the lower frequencies.
This would allow a smaller number of terms $\nnodesell$ to reach a
given accuracy for inner products.
\end{remark}              

Truncating \eqref{Jsvd} to $\nnodesell$ terms, we get,
\begin{equation}
  J_\ell(\delta k) \;\approx\; \sum_{\eta=1}^{\nnodesell} U_\eta(\delta; \ell) \Sigma_\eta(\ell) V_\eta(k; \ell)~,  \qquad \delta\in[0,\dmax],\quad k\in[0,\kmax]~,
  \label{svdtrunc}
\end{equation}
from which we construct, for each $\ell$,
\begin{equation}
Y^\svd_\eta(\vd; \ell) = U_\eta(\delta; \ell) \euler^{-\imunit \ell (\omega + \pi/2)}~, \qquad
G_\eta(k; \ell) = \Sigma_\eta(\ell) V_\eta(k; \ell)~.
\end{equation}
Note that these functions only depend on $W$.
Given these, the proposed algorithm is the following:
\begin{eqnarray}
\label{eq:fbk3}
&& \textbf{[Step 1]} \quad \fourier{Z}^\svd_\eta(q, \ell) = 2\pi \sum_{m=1}^\mmax w_m \bB(k_m; q)^* \bA(k_m; q-\ell) G_\eta(k_m; \ell) \\
&& \textbf{[Step 2]} \quad Z^\svd_\eta(\gamma; \ell) = \sum_{q=-\qmax}^\qmax \euler^{-\imunit q \gamma} \fourier{Z}^\svd_\eta(q, \ell) \\
&& \textbf{[Step 3]} \quad X^{\svd}(\vd, \gamma) = \sum_{\ell=-\ellmax}^\ellmax \sum_{\eta=1}^{\nnodesell} Y_\eta^\svd(\vd; \ell) Z^\svd_\eta(\gamma; \ell) \text{.}
\end{eqnarray}
The first step is computed directly, requiring $\bigO(\nnodes \qmax \mmax)$ operations, where $\nnodes:=\sum_{\ell}\nnodesell$ is the total number of terms.
The second step is computed with FFTs, yielding a complexity of $\bigO(\nnodes \qmax \log \qmax)$.
The third step is also computed directly with $\bigO(\ntrans \nnodes \qmax )$ operations.
Using $\qmax=\bigO(n)$ and $\mmax=\bigO(n)$,
the overall computational complexity of this scheme is therefore $\bigO(\nnodes n ^2 + \ntrans \nnodes n)$.

One key point here is the separation of $J_\ell(\delta k)$ into factors depending either on $\delta$ or on $k$.
This separation of variables allows us to efficiently calculate the sums, and the SVD accomplishes this separation with minimal error.

Another key point is that the procedure above---namely calculating a separate SVD for each $J_{\ell}(\delta k)$---is actually equivalent to a single SVD of the original translation kernel $F(\vd,\vk)$, as we now show.
Recalling the Jacobi--Anger expansion \eqref{eq:jacobi-anger},
and \eqref{Jsvd},
\begin{eqnarray}
\nonumber
F(\vd,\vk) &=& \euler^{-\imunit \vd\cdot\vk} = \euler^{-\imunit \delta k \cos(\psi-\omega)}
= \sum_{\ell\in\mathbb{Z}} J_{\ell}(\delta k) \euler^{\imunit\ell(\psi-\omega-\pi/2)} \\
&=& \sum_{\ell\in\mathbb{Z}} \euler^{\imunit\ell(\psi-\omega-\pi/2)} \sum_{\eta=1}^{\infty} U_{\eta}(\delta;\ell) \Sigma_{\eta}(\ell) V_{\eta}(k;\ell) \text{.}
\end{eqnarray}
We now relabel the singular values $\Sigma_{\eta}(\ell)$ from different
modes $\ell$ with a single new index $\zeta=1,2,\dots$,
and refer to the $\ell$ and $\eta$ corresponding to each index $\zeta$ as $\ell(\zeta)$ and $\eta(\zeta)$.
Specifically, we define
\begin{equation}
\hspace{-10mm}
U_{\zeta}(\delta) = U_{\eta(\zeta)}(\delta,\ell(\zeta)), \qquad \Sigma_{\zeta} = \Sigma_{\eta(\zeta)}(\ell(\zeta)),\qquad V_{\zeta}(k) = V_{\eta(\zeta)}(k,\ell(\zeta)) \text{.}
\end{equation}
This relabeling is chosen to have non-increasing ordering of singular values,
i.e.\
$\Sigma_{\zeta} \ge \Sigma_{\zeta+1}$ for all $\zeta=1,2,\dots$.
If $\nnodesell$ is chosen such that $\Sigma_{\nnodesell+1}(\ell) \le \epsilon$
for each $\ell$, then by construction $\Sigma_{\nnodes+1} \le \epsilon$.
Truncation at a total number of terms $\nnodes$
then involves all singular values larger than $\epsilon$,
and gives
\begin{eqnarray}
\nonumber
\hspace{-10mm}
F(\vd,\vk) \approx F^\svd(\vd, \vk) &=& \sum_{\zeta=1}^{\nnodes} \euler^{\imunit \ell(\zeta)\psi} \, \euler^{-\imunit\ell(\zeta)(\omega+\pi/2)} \, U_{\zeta}(\delta) \Sigma_{\zeta}V_{\zeta}(k) \\
&=& \sum_{\zeta=1}^{\nnodes} \euler^{-\imunit\ell(\zeta)(\omega+\pi/2)}U_{\zeta}(\delta) \, \cdot \, \Sigma_{\zeta} \, \cdot \, \euler^{\imunit\ell(\zeta)\psi}V_{\zeta}(k) \text{.}
\end{eqnarray}
The factors in the last equation correspond respectively to the left singular functions, singular values, and right singular functions of $F(\vd,\vk)$ that one would obtain by minimizing the rank-$\nnodes$ approximation
\begin{equation}
\label{svdobjective}
E_\hs^2 \; := \;
\iint_{\Omega_{\dmax}} \iint_{\Omega_{\kmax}} \biggr | F(\vd,\vk) - \sum_{\zeta=1}^{\nnodes}
U_{\zeta}(\vd) \Sigma_{\zeta} V_{\zeta}(\vk) \biggl |^{2} \mathop{d\vk d\vd} \text{.}
\end{equation}
The relationship of these vectors to those from the separate SVDs is
\begin{equation}
\hspace{-10ex}
  U_{\zeta}(\vd) = U_{\zeta}(\delta,\omega) := \euler^{-\imunit\ell(\zeta)(\omega+\pi/2)}U_{\zeta}(\delta), \qquad
V_{\zeta}(\vk) = V_{\zeta}(k,\psi) := \euler^{\imunit\ell(\zeta)\psi}V_{\zeta}(k) \text{.}
\end{equation}
So the SVD approximation described in this section is nothing more than an approximation of $F(\vd,\vk)$ built from these $U_{\zeta}$, $\Sigma_\zeta$, and $V_{\zeta}$.
As expected from the structure of the SVD, if $\zeta\neq\zeta^{\prime}$, then $U_{\zeta}(\vd)$ and $U_{\zeta^{\prime}}(\vd)$ are orthogonal over $\Omega_\dmax$, as are $V_{\zeta}(\vk)$ and $V_{\zeta^{\prime}}(\vk)$ over $\Omega_\kmax$.
Furthermore, by this orthogonality,
\begin{equation}
\label{efrob}
E_\hs^2 = \sum_{\ell \in \Integer} E_\hs(\ell)^2 = \sum_{\zeta=\nnodes+1}^\infty \Sigma_\zeta^2 \text{.}
\end{equation}

An illustration of this structure is shown in Figure~\ref{fig:Sigma}, where the array of values $\Sigma_{\eta}(\ell)$ is displayed for three choices of $W$.
Using the $W=1$ case as an example, consider what happens if we set our spectral-norm tolerance to $\epsilon=10^{-2}$ (we will make this notion more precise in the next section). 
When $\ell=0$, $\Sigma_{\eta}(0)\geq\epsilon$ only when $\eta\leq 4$; i.e., $H_{0}=4$ terms are required to achieve an $\epsilon$-accurate approximation to $J_{0}(\delta k)$.
On the other hand, when $\ell=3$, $\Sigma_{\eta}(3)\leq \epsilon$ only when $\eta\leq 2$; i.e., only $H_{3}=2$ terms are required for an $\epsilon$-accurate approximation to $J_{3}(\delta k)$.
In order to combine the SVD expansions of the various $J_{\ell}(\delta k)$ into an $\epsilon$-accurate expansion of $F(\vd,\vk)$, we need to retain all the terms $\Sigma_{\eta}(\ell)\geq\epsilon$.
With this choice of $\epsilon$, we retain the $\nnodes=34$ largest values of $\Sigma_{\zeta}$, which correspond to the $\nnodes = \sum_\ell \nnodesell = 34$ largest values of $\Sigma_{\eta}(\ell)$.
In the next section will use this structure to discuss the approximation error of this expansion.

\begin{figure}
\centering
  \includegraphics[width=0.9\linewidth,trim= {0cm 0cm 0cm 0cm},clip]{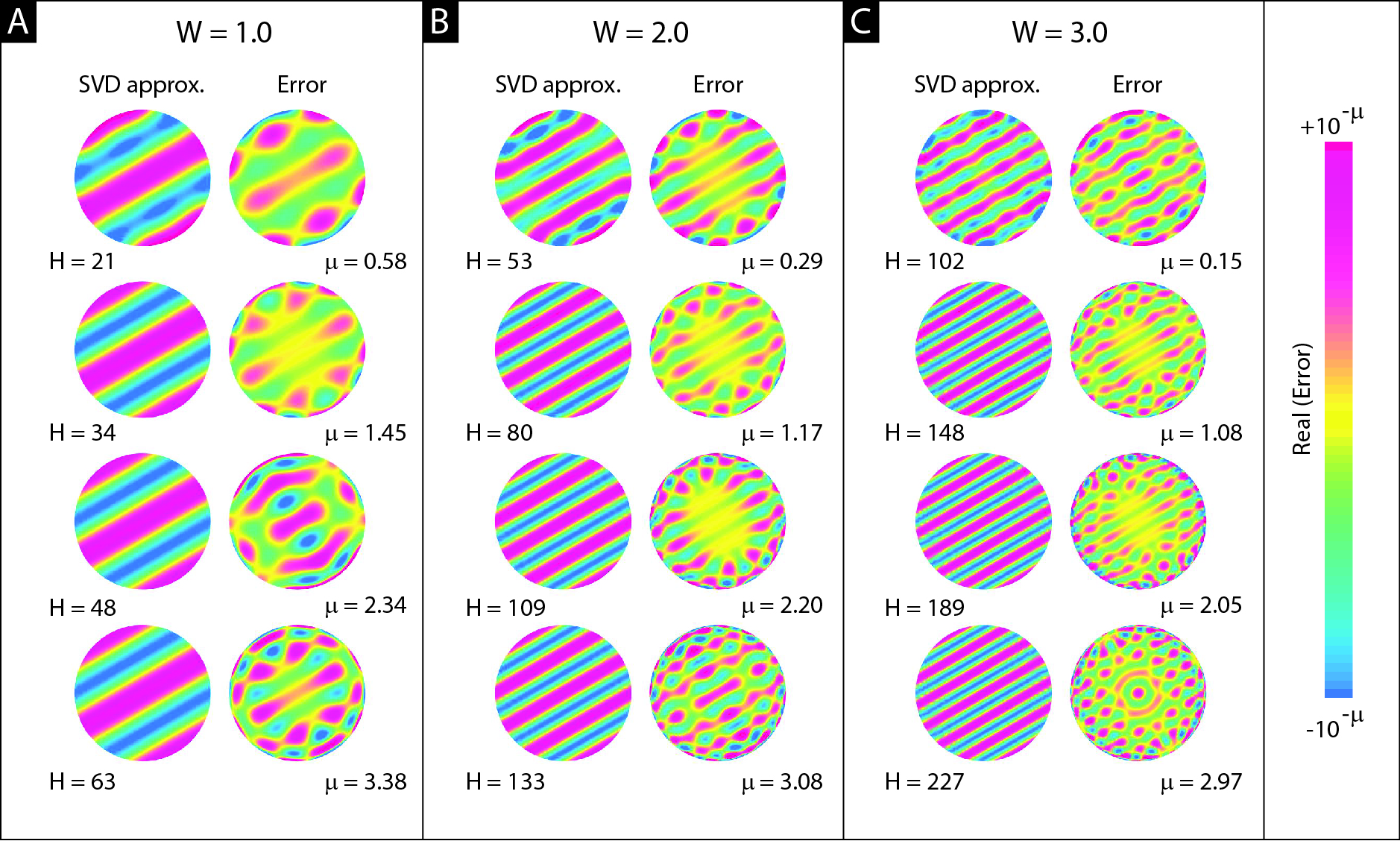}
\caption{
\label{fig:xn}
Here we illustrate the SVD approximation to the translation kernel $F\left(\vd,\vk \right) =\euler^{-\imunit\vd\cdot\vk}$, with $\vd$ chosen so that $\omega=\pi/4$ and $\delta=\dmax = 2\pi \wmax/\kmax$.
In panel A we show a pair of columns associated with $W=1$.
Each row in this panel shows the real component of the SVD approximation on the left and the real component of the associated error on the right.
The number of terms $\nnodes$ increases with each row.
Both the left and right columns use the same color scale (far right), but with different limits:\ the SVD approximation to $F$ has limits of $[-1, 1]$, whereas the error has limits of $[-10^{-\mu}, 10^{-\mu}]$, with the value of $\mu $ listed below each row.
The value of $\mu$ is the number of digits of (pointwise) accuracy achieved by the corresponding SVD approximation.
Panel B and panel C display similar results for $W=2$ and $3$.
}
\end{figure}

We illustrate this SVD approximation in Figure~\ref{fig:bessel-svd} which shows the structure of $U_{\eta}(\delta; 0)$ and $V_{\eta}(k; 0)$ for the lowest-order Bessel function $J_0(\delta k)$.
Figure~\ref{fig:xn} illustrates the errors incurred by the SVD expansion when approximating the translation kernel $F(\vd, \vk)$.
Note that the SVD approximation provides a successively more accurate approximation of $F$ as the number of terms $\nnodes$ increases.
In addition, as $W$ increases, a larger $\nnodes$ is required to achieve the same level of accuracy (compare Figure~\ref{fig:xn}A with Figure~\ref{fig:xn}C).
The next section places rigorous bounds on this growth.

\section{Error analysis\label{sec:erroranalysis}}

In this section we prove a rigorous upper bound on
the number of singular values needed in \meth
to achieve a relative error $\epsilon$ in inner product computations,
given a maximum dimensionless translation $W$.
We first prove the bound on the $\epsilon$-rank of the translation
kernel, then show how this controls
the relative error over the set of inner products.

\subsection{A bound on the $\epsilon$-rank of $F(\vd, \vk)$}

Here we prove that the number of retained SVD terms
need grow at worst quadratically
in $W$ and linearly in the desired number of digits of accuracy $\log(1/\epsilon)$.

\begin{thm}  
  Let $W = \dmax\kmax / 2\pi$ be the maximum translation
  in wavelengths defined in Section~\ref{sec:setup},
  and let $\epsilon>0$ be the desired precision.
  Let $H$ be the number of singular values $\Sigma_\zeta$
  of the translation kernel $F(\vd,\vk)$ larger
  than $\epsilon$, i.e.\ $H$ is the $\epsilon$-rank of $F$.
  Then
  \begin{equation}
    H = \bigO\bigl((W + \log(1/\epsilon))^2\bigr)~.
  \end{equation}
  \label{t:rank}
\end{thm}
    {\bf Proof.}  
    As discussed in the previous section,
    if $\nnodesell$ is the number of singular values
    $\Sigma_{\eta}(\ell)$ greater than $\epsilon$, for the
    translation kernel with kernel $J_\ell(\delta k)$,
    then 
    $\nnodes$ in the theorem statement equals $\sum_\ell \nnodesell$.
    The lemma below implies that this is achieved with
    $0\le \nnodesell \le \max(0, {\cal H}_{\epsilon,W} - |\ell|/2)$, with
    ${\cal H}_{\epsilon,W} := \max(\pi \euler^2 W, \log(2\pi W/\epsilon)) +3/2$.
    This means that $H_\ell = 0$ for $|\ell| > 2{\cal H}_{\epsilon,W}$.
    Then $H = \sum_{\ell \in \Integer} \nnodesell \le \sum_{|\ell|\le 2{\cal H}_{\epsilon,W}} ({\cal H}_{\epsilon,W} - |\ell|/2) = 2{\cal H}_{\epsilon,W}^2$.
    Asymptotically, we may summarize
    ${\cal H}_{\epsilon,W} = \bigO(W + \log(1/\epsilon))$, completing the proof.
\hfill$\square$  

The key lemma below shows that the number of terms $\nnodesell$
grows linearly with $W$ and $\log(1/\epsilon)$,
but shrinks linearly with $|\ell|$.
This behavior is clear in Figure~\ref{fig:Sigma}
as the roughly triangular shape of the set of indices with
significant singular values.

\begin{lem}    
  Let $\ell\in\mathbb{Z}$,
  and let $\nnodesell \ge 0$ be the $\epsilon$-rank of the modal
  translation kernel
  $J_\ell(\delta k)$ acting from $L^2([0,\kmax])$ to
  $L^2([0,\dmax])$ with
  linear weight functions as defined in the previous section.
  Then $\nnodesell \; \le \; \overline{H}_\ell$, where
  \begin{equation}
    \overline{H}_\ell
    \; := \;
    \max\bigl(0, \, \pi \euler^2 W - |\ell|/2, \, \log(2\pi W/\epsilon) + 3/2 - |\ell|/2 \bigr)~.
    \label{rank}
  \end{equation}
\label{l:rank}
\end{lem}
{\bf Proof.}       
Defining rescaled variables $x := \kmax \delta$ and $y:=k/\kmax$,
the kernel becomes $J_\ell(\delta k) = J_\ell(x y)$ from functions
over $y \in [0,1]$ to those over $x\in [0, 2\pi W]$;
one may check that rescaling does not change the singular values $\Sigma_\eta(\ell)$.
Consider the case $\ell$ even.
We exploit an identity given by Wimp \cite[Eq.~(2.22)]{wimp},
\begin{equation}
  \hspace{-10ex}
  J_\ell(xy) = \sum_{\nu=0}^\infty C_{2\nu}(x) T_{2\nu}(y), \quad |y|\le 1,
  \qquad C_{2\nu}(x)= \varepsilon_{\nu} J_{\ell/2 + \nu}(x/2) J_{\ell/2 - \nu}(x/2),
  \label{wimp}
  \end{equation}
where $T_{2\nu}$ is the $2\nu$th Chebyshev polynomial,
and $\varepsilon_{\nu} = 1$ for $\nu=0$ and 2 otherwise.
Applying $|J_{|\ell|/2-\nu}(x/2)|\le 1$ \cite[Eq.~(10.14.1)]{dlmf} to the
smaller magnitude of the two Bessel orders,
we get a bound in terms of the Bessel with larger magnitude order
\begin{equation}
  |C_{2\nu}(x)| \le 2 | J_{|\ell|/2+\nu}(x/2) |~.
\label{C2nx}
\end{equation}

Recalling the definition in \eqref{rank}, we now claim that
\begin{equation}
  \sum_{\nu\ge\overline{H}_\ell} |C_{2\nu}(x)| \le \frac{\epsilon}{2\pi W}~,
  \qquad \mbox{ for all } x\in[0,2\pi W]~.
  \label{claim}
\end{equation}
To prove this, for each $\nu$ we rescale the Bessel argument in \eqref{C2nx}
by writing $p := x/(|\ell|+2\nu)$
as the ratio of argument $x/2$ to order $|\ell|/2+\nu$,
or fraction of the distance from the origin to the Bessel turning point.
We then apply Siegel's
bound \cite[Eq.~(10.14.5)]{dlmf}
in the evanescent region $0\le p \le 1$,
\begin{equation}
  \hspace{-10ex}
|J_{|\ell|/2+\nu}((|\ell|/2+\nu)p)| \;\le\;
\euler^{(|\ell|/2 + \nu) \,\left(\log p + \sqrt{1-p^2} - \log\left(1+\sqrt{1-p^2}\right)\right)}~,
\label{siegel}
\end{equation}
From \eqref{rank}, we obtain that $\nu \ge \overline{H}_\ell$ implies $2\pi W \le \euler^{-2} (|\ell| + 2\nu)$.
Since $x \le 2\pi W$, we then have $p \le \euler^{-2}$ and $0.99 < \sqrt{1-p^2} < 1$.
The second factor in the exponent of \eqref{siegel} is therefore less than $-2 + 1 - \log(1.99) < -1$.
Consequently, \eqref{siegel} is bounded by simple exponential decay:
\begin{equation}
  |J_{|\ell|/2+\nu}(x/2)| \;\le\; \euler ^{-(|\ell|/2 + \nu)}~,
  \quad \mbox{ for all } x\in[0,2\pi W], \quad \nu  \ge \overline{H}_\ell.
\label{decay}
\end{equation}
This bounds the left side of \eqref{claim} by the geometric series,
$$
\sum_{\nu\ge\overline{H}_\ell} |C_{2\nu}(x)| \;\le\;
2 \sum_{\nu\ge\overline{H}_\ell} \euler ^{-(|\ell|/2 + \nu)} \;\le\;
\frac{2}{1-\euler^{-1}} \euler^{-(|\ell|/2 + \overline{H}_\ell)}~.
$$
Substituting the bound
$\overline{H}_\ell \ge \log(2\pi W/\epsilon) + 3/2 - |\ell|/2$
from \eqref{rank}, and using $\euler^{3/2} > 2/(1-\euler^{-1})$,
now establishes the claim \eqref{claim}.

We work in weighted $L^2$-norms,
$\|f\|^2 := \int_0^a |f(x)|^2 \mathop{x dx}$, where $a=1$
for the domain or $a=2\pi W$ for the range.
By Allahverdiev's theorem \cite[p.~28]{GK69},
the error in the induced operator norm $\| \cdot \|$
of any rank-$r$ approximation is at least the $(r+1)$th
operator singular value.
Choosing $r=\overline{H}_\ell$ and recognizing
the $\nu<r$ truncation of \eqref{wimp} as a particular rank-$r$ approximation
to $J_\ell$, we insert the definition of the operator norm and get
\begin{eqnarray}
\hspace{-10mm}
\Sigma_{\overline{H}_\ell+1}(\ell)^2 &\le\;
\biggl\| J_\ell(\cdot,\cdot) -
\sum_{\nu=0}^{\overline{H}_\ell-1} C_{2\nu}(\cdot) T_{2\nu}(\cdot) \biggr\|^2
\; = \;
\biggl\| \sum_{\nu\ge\overline{H}_\ell} C_{2\nu}(\cdot) T_{2\nu}(\cdot) \biggr\|^2
\nonumber
\\
&= \sup_{\|f\|=1} \int_0^{2\pi W} \biggl|
\int_0^1 \sum_{\nu \ge \overline{H}_\ell} C_{2\nu}(x) T_{2\nu}(y)f(y) \mathop{y dy} \biggr|^2 \mathop{x dx}
\nonumber
\\
& \le \;
\sup_{\|f\|=1} \int_0^{2\pi W}
\biggl(\int_0^1 \biggl| \sum_{\nu \ge \overline{H}_\ell} C_{2\nu}(x) T_{2\nu}(y)\biggr|^2 \mathop{y dy} \biggr) \cdot
\biggl(\int_0^1 f(y)^2 y dy\biggr) \mathop{x dx}
\nonumber \\
& \le \;
\int_0^{2\pi W}
\biggl(\sum_{\nu \ge \overline{H}_\ell} |C_{2\nu}(x)|\biggr)^2 \cdot 1 \, \mathop{x dx}
\; \le \;
(2\pi W)^2 \biggl(\frac{\epsilon}{2\pi W}\biggr)^2 \; \le \; \epsilon^2~,
\end{eqnarray}
where to reach the third line we used the $y$-weighted Cauchy--Schwarz inequality,
to reach the fourth line we used the fact that the magnitude of the Chebyshev polynomials are bounded by one,
and finally the claim \eqref{claim}.
Note that, in second expression above, when $\nnodesell=0$ there are no terms in
the sum.
The bound just shown, $\Sigma_{\overline{H}_\ell+1}(\ell) \le \epsilon$,
is equivalent to the statement
that the $\epsilon$-rank is no more than $\overline{H}_\ell$,
completing the proof.

The proof for $\ell$ odd is similar, exploiting the identity \cite[Eq.~(2.23)]{wimp}.
\hfill$\square$   

\begin{remark}
  The above proof gives explicit constants in the upper bound on $\nnodesell$,
  but since a small argument-to-order ratio $p<\euler^{-2}$
  was chosen, these constants are far from being tight and could be improved.
  The empirical behavior of $\nnodesell$ has the same linear
  form as \eqref{rank} but with different constants: a numerical fit gives
  $
  \nnodesell \approx 2.0 (\wmax - |\ell|/2\pi) + 0.5 \log(1/\epsilon)
  $.
  This is strong evidence that the quadratic power in Theorem~\ref{t:rank} is tight.
\end{remark}

\subsection{Induced error on inner product $X(\vd, \gamma)$}
\label{sec:Xbound}

Consider the induced, or spectral, norm $E_2$ of the difference between the
translation kernel and its rank-$\nnodes$ expansion, defined by
\begin{equation}
  E_2^2 := \sup_{\|\fourier{A}\|_2 = 1} \iint_{\Omega_\dmax} \left | \iint_{\Omega_\kmax}
  \biggl(F(\vd,\vk) - F^\svd(\vd,\vk) \biggr)
  \fourier{A}(\vk) \mathop{d\vk} \right |^2 \mathop{d\vd} \text{,}
\end{equation}
where $\fourier{A} \in \Leb^2(\Omega_\kmax)$, and $\|\fourier{A}\|_{2}$ is the 2-norm on $\Leb^2(\Omega_\kmax)$.
By Allahverdiev's theorem \cite[p.~28]{GK69}
(see the proof of Lemma~\ref{l:rank}),
\begin{equation}
  E_2 = \Sigma_{\nnodes+1}~.
  \end{equation}
We now show that \meth induces a small error on the inner product $X(\vd, \gamma)$.
To simplify the analysis, we consider the \meth approximation of its continuous analog ${\cal X}(\vd, \gamma)$,
\begin{equation}
    {\cal X}^\svd(\vd,\gamma) = \iint_{\Omega_\kmax} F^\svd(\vd,\vk) \fA(\vk)^\ast R_{-\gamma} \fB(\vk) \mathop{d\vk} \text{.}
\end{equation}
The squared $\Leb^2$-norm of the residual ${\cal X}(\vd, \gamma) - {\cal X}^\svd(\vd, \gamma)$ as a function over $\Omega_\dmax \times [0, 2\pi)$ is now bounded as follows:
  \begin{eqnarray}
    \hspace{-12ex}
    \left \| {\cal X}-{\cal X}^\svd \right\|_2^2 &=& \int_{0}^{2\pi} \iint_{\Omega_\dmax} \left| {\cal X}(\vd, \gamma)-{\cal X}^\svd(\vd, \gamma) \right|^2 \mathop{d\vd} \mathop{d\gamma}
    \nonumber \\
&=& \int_0^{2\pi} \iint_{\Omega_\dmax} \left| \iint_{\Omega_\kmax}
\biggl(F(\vd,\vk) - F^\svd(\vd,\vk)  \biggr)
\fA(\vk)^\ast R_{-\gamma} \fB(\vk) \mathop{d\vk} \right|^2 \mathop{d\vd} \mathop{d\gamma} \nonumber \\
&\leq& E_2^2 \int_0^{2\pi} \iint_{\Omega_\kmax} \left| \fA(\vk)^\ast R_{-\gamma}\fB(\vk) \right|^2 \mathop{d\vk} \mathop{d\gamma} \nonumber \\
&=& \Sigma_{\nnodes+1}^{2} \int_0^{2\pi} \| \fA(\cdot)^{\ast} R_{-\gamma}\fB(\cdot) \|_{2}^{2} \mathop{d\gamma} \text{.}
\label{Xmserr}
  \end{eqnarray}
  Thus, the relative root mean square error in the inner product function
  is controlled by $\Sigma_{H+1}$,
  which, if $H$ is chosen as in Theorem~\ref{t:rank}, cannot
  exceed the desired $\epsilon$.
  The last factor in \eqref{Xmserr}
  depends only on the overlap of the radial power spectra of $A$ and $B$ and may be rewritten as $2\pi \int_0^K |\bA(k, 0) \bB(k, 0)|^2 \mathop{dk}$.
Note that the bound \eqref{Xmserr} is tight, and is achieved if $\fA(\vk)^\ast R_{-\gamma}\fB(\vk)$ equals the first omitted right singular function $V_{\nnodes+1}(\vk)$.

The above argument provides a bound on the error induced by using the truncated translation kernel $F^\svd(\vd, \vk)$ in ${\cal X}(\vd, \gamma)$.
The discretized inner products $X(\vd, \gamma)$ and $X^\svd(\vd, \gamma)$ differ from ${\cal X}(\vd, \gamma)$ and ${\cal X}^\svd(\vd, \gamma)$ only by quadrature and Fourier truncation errors, which can be made arbitrarily small (see Section~\ref{sec:discretization}).
Consequently, the error $\|X(\vd, \gamma) - X^\svd(\vd, \gamma)\|_2$ is dominated by $\Sigma_{H+1}$, and hence $\epsilon$.

\section{Numerical results}
\label{sec:results}

\begin{figure}
\centering
\raisebox{-.2in}{
  \includegraphics[width=.25\linewidth,trim= {2.5cm 1.5cm 2.5cm 1.0cm},clip]{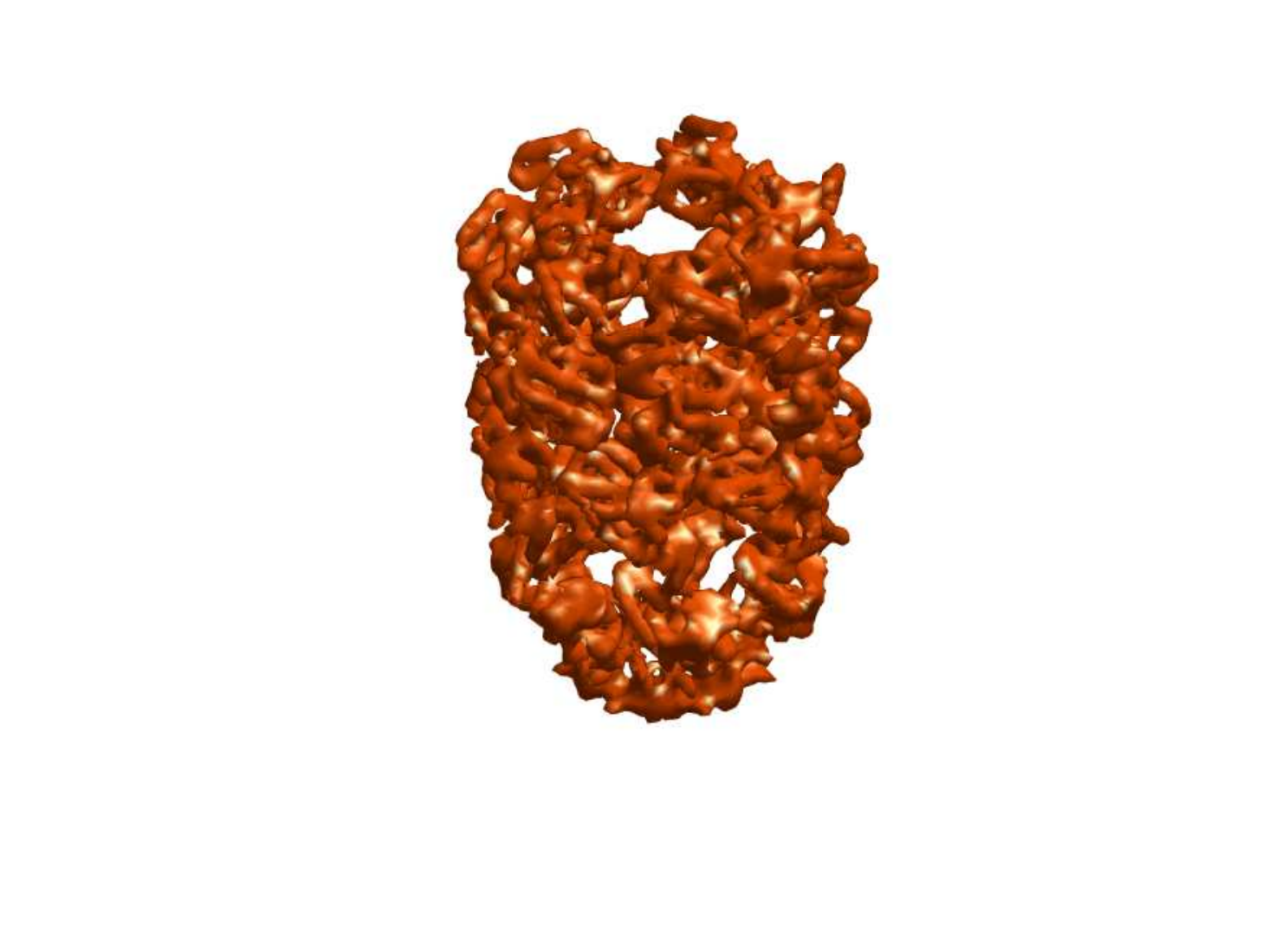}}%
\hspace{2mm}%
\fbox{\includegraphics[width=.15\linewidth,trim= {1.0cm 1.0cm 1.0cm 1.0cm},clip]{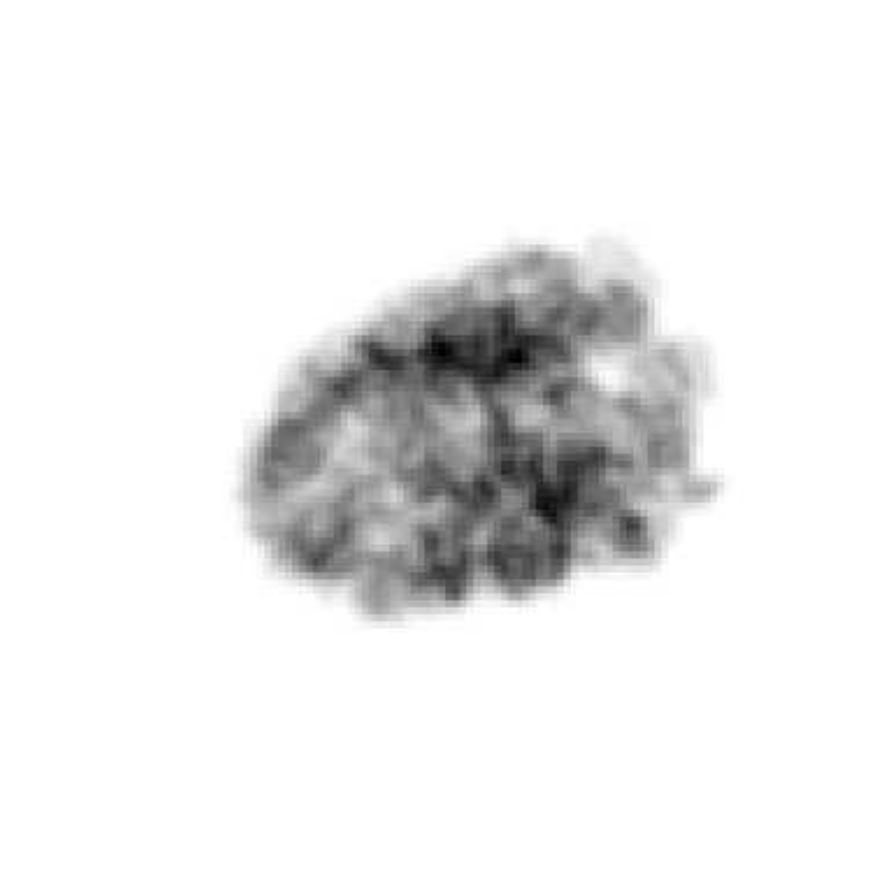}}%
\hspace{2mm}%
\fbox{\includegraphics[width=.15\linewidth,trim= {1.0cm 1.0cm 1.0cm 1.0cm},clip]{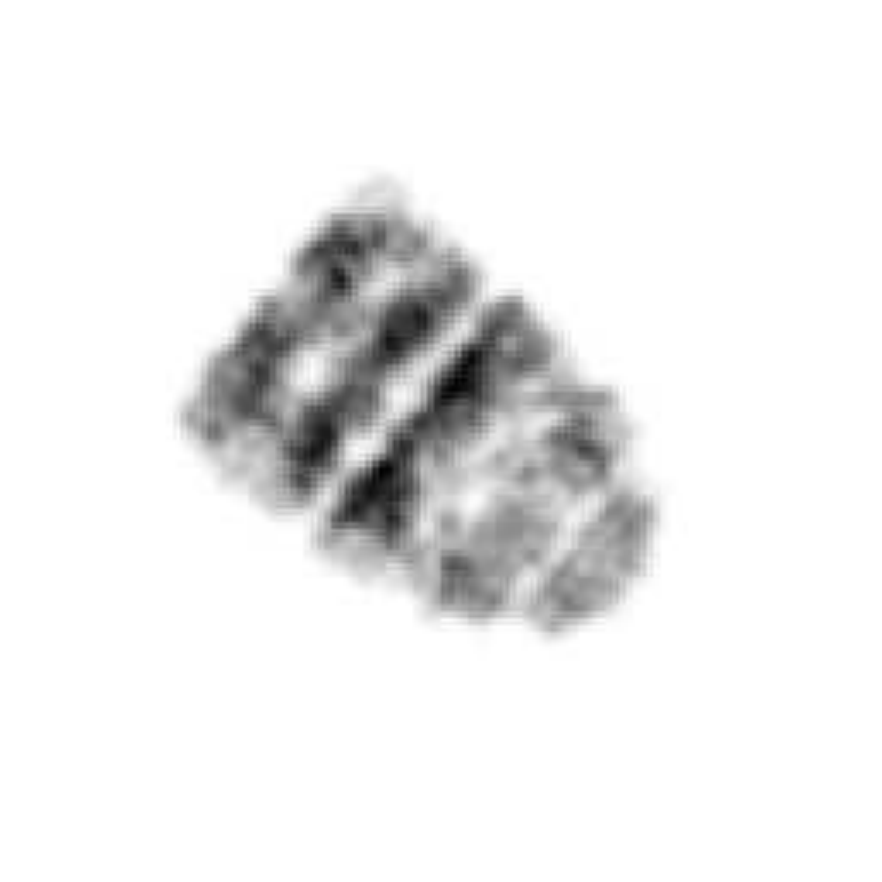}}%
\hspace{2mm}%
\fbox{\includegraphics[width=.15\linewidth,trim= {1.0cm 1.0cm 1.0cm 1.0cm},clip]{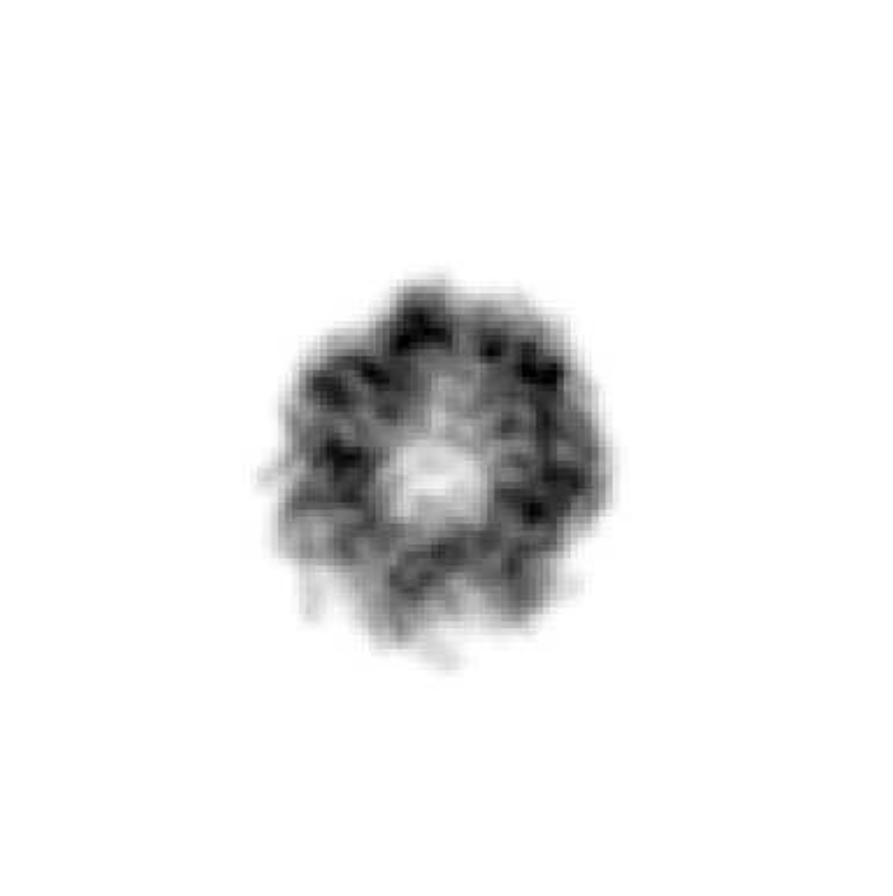}}%
\hspace{2mm}%
\fbox{\includegraphics[width=.15\linewidth,trim= {1.0cm 1.0cm 1.0cm 1.0cm},clip]{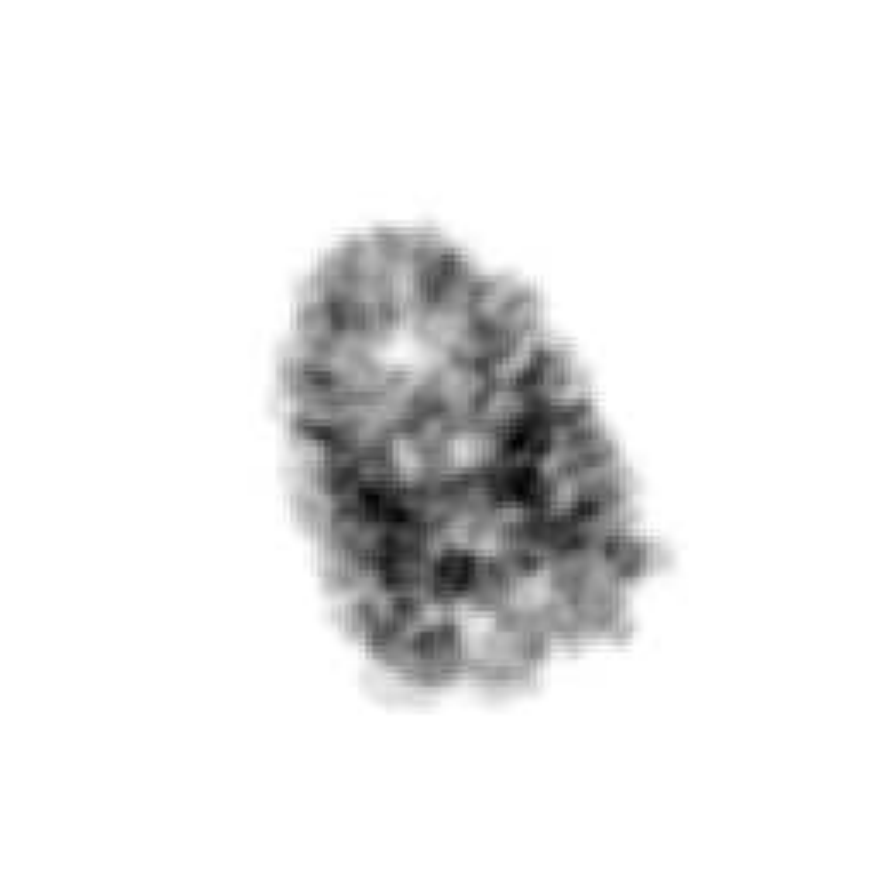}}%
\caption{Four projection templates of the 3D electron density shown on the
  far left, taken at random Euler angles, and sampled on a $128\times 128$ pixel grid.
  The images are then obtained
  by randomly rotating and translating these 2D templates
  with sub-pixel resolution.
    The task is to recover the true translations and rotations as well as the templates that match each of the images.
\label{figure:templates}}
\end{figure}
In this section we compare the proposed \meth method of Section~\ref{sec:svd} with two existing fast methods for evaluating the inner product function
$X(\vd,\gamma)$: the BFR and BFT methods (see Section~\ref{sec:intro}).
The methods are compared for the task of rigidly aligning $\nim \gg 1$ images to another set of $\nim$ images, referred to as ``templates'' for clarity.
We first define the BFR and BFT methods and describe their implementations.

\subsection{Brute force rotations (BFR)}
\label{sec:bfr}

In short, for each of the $\nim^2$ image--template pairs, and each of the $\bigO(n)$ rotations of one of the images, we perform a standard fast 2D convolution of $A$ with $R_{-\gamma} B$ using 2D FFTs.
This relies on writing \eqref{Xdg} as follows, using the 2D convolution theorem \cite{bracewell}:
\begin{equation}
\hspace{-0.75in}
{\cal X}(\vd,\gamma)=
\iint_{\Real^2} A(\vx-\vd) R_{-\gamma} B(\vx)^\ast d\vx
= \frac{1}{(2\pi)^2} \iint_{\Omega_\kmax} \fourier{A}(\vk)\, R_{-\gamma}\fourier{B}(\vk)^\ast
\, \euler^{-\imunit \vk \cdot \vd} d\vk~.
\label{2dconv}
\end{equation}
Discretizing the right-hand side on a Cartesian grid, we may now recover ${\cal X}(\vd,\gamma)$ on an entire $n\times n$ Cartesian $\vd$ grid with a single 2D FFT.
Since a 2D FFT is needed for each of $\bigO(n)$ rotations and $\nim^2$ image--template pairs, the cost per pair is $\bigO(n^3\log n)$, which is one power of $n$ slower than the proposed SVD method; see Table~\ref{table:complexity}.
A finer sampling in $\vd$ is achieved by zero-padding each FFT.
Note that the cost is independent of $\wmax$ or $\ntrans$.

For each of the $\nim$ images $A$ and for each of the $\bigO(\nim n)$ rotated templates $R_{-\gamma} B$, a precomputation is needed to evaluate $\fourier{A}(\vk)$ and $R_{-\gamma}\fourier{B}(\vk)$ for each $\gamma$ on an $n\times n$ Cartesian grid in $\vk$.
A common approach for computing $R_{-\gamma}\fourier{B}(\vk)$ on such a grid is to resample by interpolation in real or Fourier space.
We instead use the 2D type 2 NUFFT \cite{finufft} as in Section~\ref{sec:nufft} to get Fourier transform samples on polar grids with spectral accuracy;
this preparation step is not included in our timings.
Then for each $\gamma$ we rotate on the polar grid and
use a 2D type 1 NUFFT to resample to a Cartesian Fourier grid.
We do include the latter in timings, but since this precomputation scales
only as $\nim$, it is practically negligible compared to the main computation.

\subsection{Brute force translations (BFT)}
\label{sec:bft}

For each of the $\nim$ templates $B$ we
precompute their Fourier--Bessel coefficients $\bB(k;q)$
on the rings $\{k_m\}_{m=1}^\mmax$, $|q|\le \qmax$ (see Section~\ref{sec:nufft}).
For each of the $\nim$ images $A$ and for each of their $\ntrans$ desired
translations, we use the same method for the
translated Fourier--Bessel coefficients $T_{\vd}\bA(k;q)$,
except we multiply $\fA$ by the translation phase \eqref{Fdk}
before the 1D FFTs in \eqref{aqkcalc}.
This precomputation requires $\bigO(\nim \ntrans n^2\log n)$ operations.

For each of the $\nim^2 \ntrans$ pairs
of templates and $\vd$-translated images,
we then apply \eqref{eq:fbk2} with $\ellmax=0$ to calculate
$\hat X(\vd,q)$ and use \eqref{1dfft} as in Section~\ref{sec:fastrots}.
This returns all inner products associated with $\vd$ for the grid of
$\gamma\in[0,2\pi)$.
The total cost of precomputation is
$\bigO(\nim^2 \ntrans (n^2 + n\log n))$ effort, i.e.\ $\bigO(\ntrans n^2)$
per image--template pair.
We remark that a similar method was used for efficiently
computing inner products across multiple rotations in \cite{Barnett2017}.
Note that Joyeux {\it et al.} \cite{Joyeux2002} instead use
    angular convolutions in real space, but that the expected cost is similar.

\subsection{Performance comparison for cryo-EM template matching}

To conduct our numerical experiments, we start with the GroEL--GroES 3D structure \cite{Xu1997}.
We compute $\nim=10$ random tomographic projections (see Figure~\ref{figure:templates}) to use as templates.
The size of each template is $128 \times 128$ pixels (i.e., $n = 128$), and the templates have unit $L^2$-norm in $[-1,1]^2$.
The maximum frequency $\kmax$ is chosen as the Nyquist frequency (see Remark~\ref{r:nyq}).
To generate images, we rotate the templates by random real angles in $[0, 2\pi)$ and translate them by random shifts within $\Omega_\dmax$, where $\dmax$
is varied in a manner described below.
  All such projections and transformations
  are computed using Fourier methods which
  perform spectrally accurate off-grid translations; see \cite{Barnett2017}.
  
  The timing results below are carried out with an efficient vectorized MATLAB
  implementation.
  Specifically,
  for BFT and \meth, computations are dominated by inner products,
  which are blocked together as matrix-matrix multiplications.
  Since these exploit level 3 BLAS, they are essentially as efficient as possible
  on the architecture.
  For BFR, the cost is dominated by FFT calls, which
  use MATLAB's version of the FFTW library, and are thus
  also very efficient.
  We run MATLAB 2016b on a Linux workstation with two six-core Intel Xeon E5-2643 CPUs
  at 3.4 GHz and 128 GB of memory.
  
We perform experiments over a wide range of $\dmax$.
For a given $\dmax$, which in turn fixes $\wmax = n\dmax/4$, we cover $\Omega_\dmax$ with an approximately uniform grid of translations $\vd$ of two different realistic mean densities:
\begin{enumerate}
\item[(a)] a coarse grid of 4 translations per square pixel, i.e.\ mean spacing $\dx/2 \approx 0.0078$, yielding $N \approx \pi n^2 \dmax^2$ translations, and
\item[(b)] a finer grid of 16 translations per square pixel, i.e.\ mean spacing $\dx/4 \approx 0.0039$, yielding $N \approx 4\pi n^2 \dmax^2$ translations.
\end{enumerate}
The former is commonly used in applications, while in our setting the latter is
sufficient for robust alignment to the desired accuracy.
We also search over a rotation grid of $2\pi \kmax = n \pi^2 \approx 1300$
equispaced $\gamma$ values; this grid does not depend on $\dmax$.
We now loop through image--template pairs, computing the
array of inner products $X(\vd,\gamma)$ on the translation and
rotation grids using the three methods described above: BFR, BFT, and
the proposed \meth method.
For each experiment, we record the full time taken to compute these inner product arrays. 
The full computation time includes the precomputation scaling as $\nim$, which involves preparing
the images and templates for inner product calculations, plus the actual computation of the inner products scaling as $\nim^2$.
The latter dominates the cost as $\nim$ grows.
We exclude from our timing results the one-time \meth planning stage that is independent of the images, such as filling the translation kernel matrices and computing their SVD; this is negligible in any realistic application.
In \meth, we set the (relative) accuracy to $\epsilon = 10^{-2}$, which, together with $\wmax$, determines the number of SVD terms $\nnodes$.

As a measure of consistency, we ensure that---after the calculation---the maximal inner product for each image corresponds to that of the true template and the nearest possible on-grid shift and angle used to generate that image.

\begin{figure}
\centering
  a)\raisebox{-2.2in}{
    \includegraphics[width=0.45\linewidth,trim= {0cm 0cm 0cm 0cm},clip]{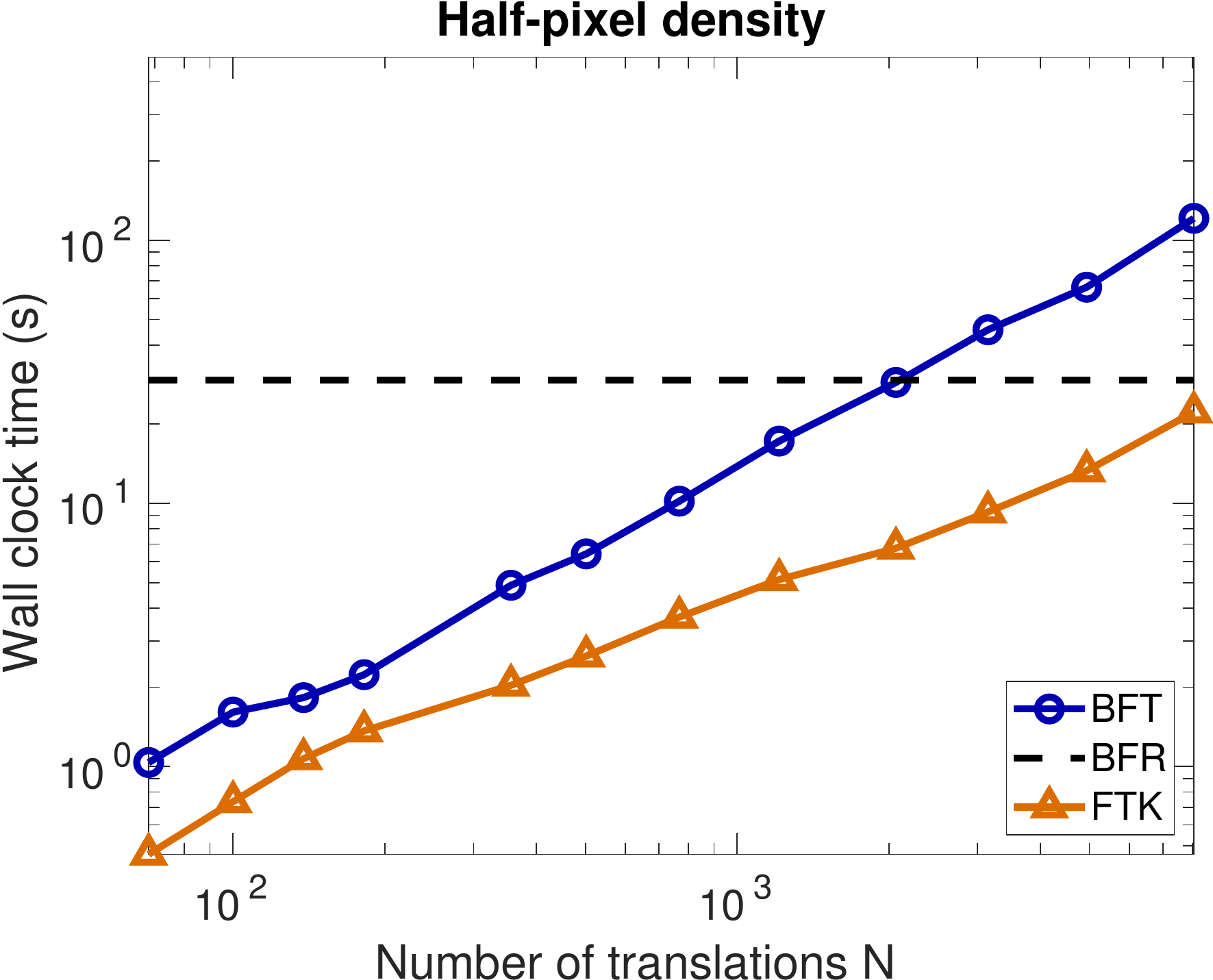}
    }
  b)\raisebox{-2.2in}{
    \includegraphics[width=0.45\linewidth,trim= {0cm 0cm 0cm 0cm},clip]{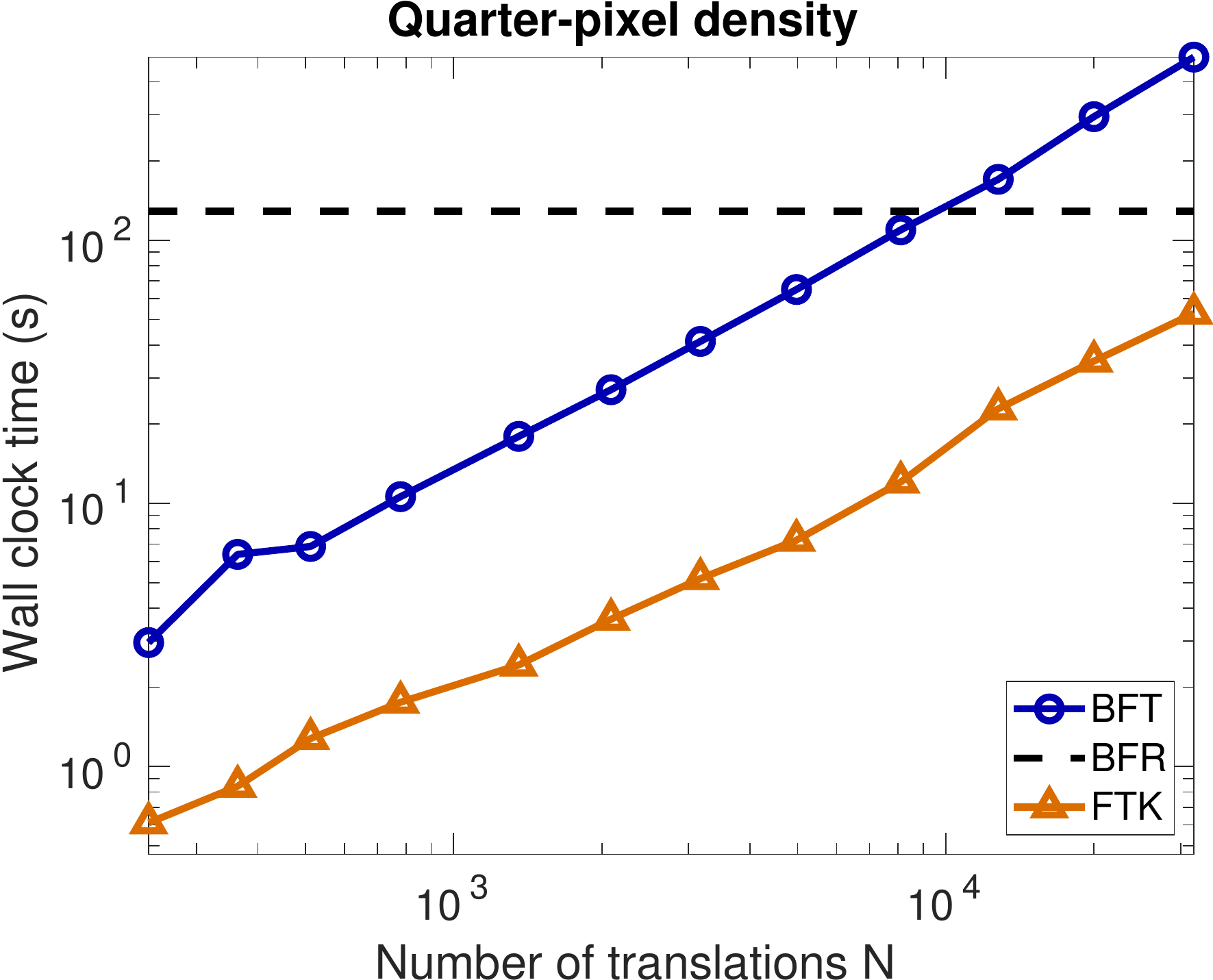}
    }
\caption{ \label{figure:results} 
Timing results for the calculation of
  all inner products between rigid transformations of $\nim=10$ images
  and $\nim$ templates, of size $n\times n$ pixels for $n=128$.  Three
  methods are compared: BFR (brute force rotations, as
  implemented in Section~\ref{sec:bfr}), BFT (brute force
  translations, as in Section~\ref{sec:bft}), and our
  proposed \meth method (Section~\ref{sec:svd}).
  The horizontal axis shows $\ntrans$; this scales as
  $\bigO(\dmax^2)$ because we choose a constant density of
  translations of in a disk of radius $\dmax$.  The mean translation spacings are (a) half-pixel $\dx/2$ and (b) quarter-pixel $\dx/4$.
  In both plots, the range of $\dmax$ is $[0.04, 0.4]$, i.e., between $1/50$ and $1/5$ of the image size.  }
\end{figure}

Figure~\ref{figure:results} compares the running time of the BFR and
BFT methods with our proposed \meth method across a range of $\ntrans$
corresponding to $\dmax$ between $0.04$ and $0.4$.
The largest $\dmax$ is equivalent to a $25$-pixel shift.
The comparison shows that, for a wide range of
$\ntrans$, \meth outperforms BFT by a factor of $3$ for the $\dx/2$ spacing.
For the $\dx/4$ spacing grid, we observe a speedup factor of $8$--$10$.
At large $\ntrans$, the precomputation stage of BFT takes about one third of the overall running time, while for \meth, roughly one tenth of the time is spent in precomputation.
The figure also shows that BFR takes about $30~\mathrm{s}$ for $\dx/2$ and $100~\mathrm{s}$ for $\dx/4$, independent of $\ntrans$.
In both cases, for all $\dmax \le 0.4$ (shifts up to $25$ pixels),
\meth is faster than BFR.

\begin{figure}
\centering
  \includegraphics[width=0.6\linewidth,trim= {0cm 0cm 0cm 0cm},clip]{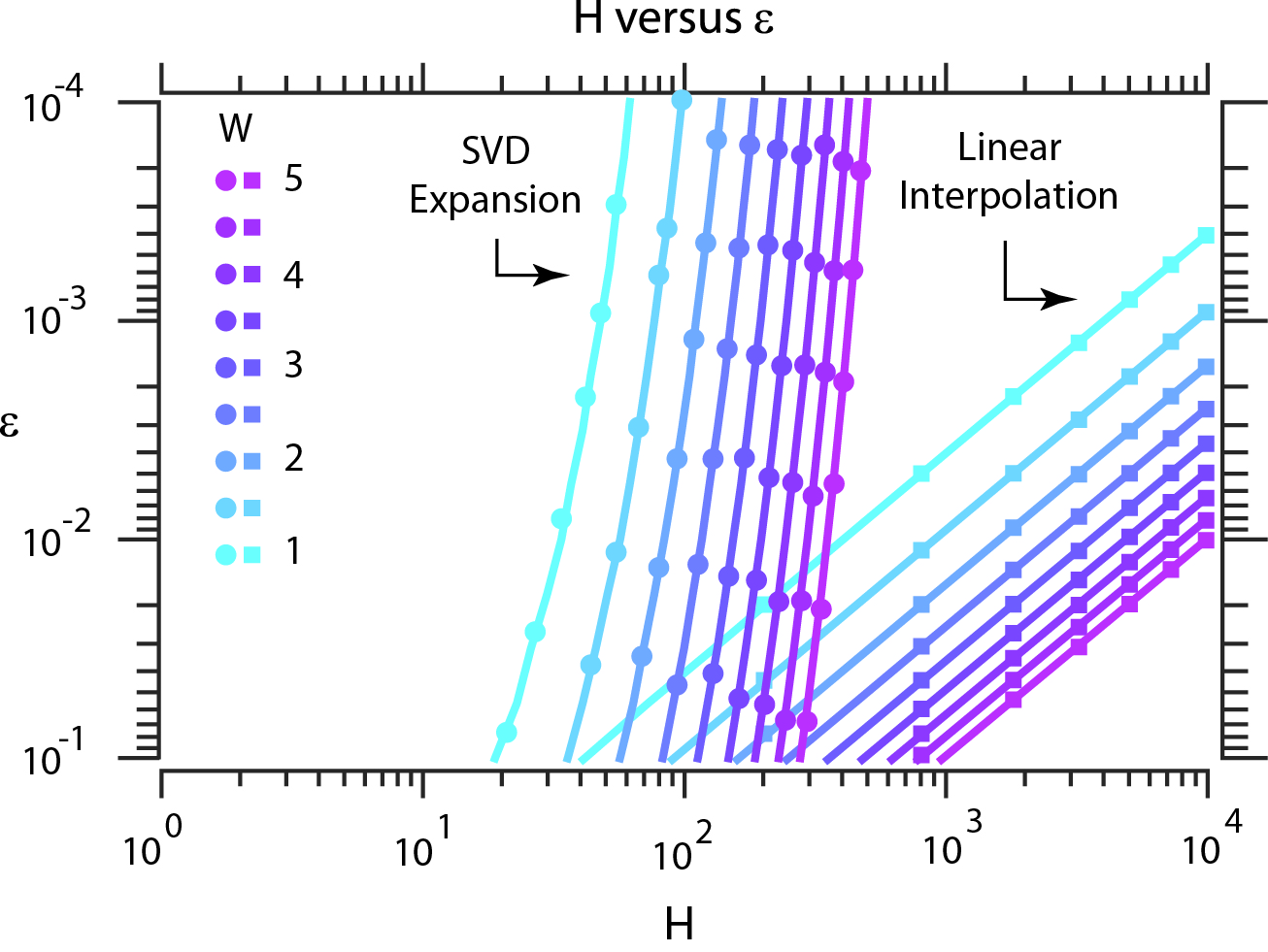}
\caption{ \label{fig:Hvseps}
The accuracy of our SVD expansion (solid dots) as a function of the number of SVD terms $\nnodes$, on logarithmic scales.
We measure accuracy in terms of $\varepsilon = E_\hs$ of \eqref{svdobjective}.
By \eqref{efrob} and Section~\ref{sec:Xbound}, this provides a bound on $\epsilon$.
The value $\log_{10}(1/\varepsilon)$ can be interpreted as the number of digits of accuracy of our expansion.
Results for different values of $W$ are shown in different colors (see legend at the left).
Results for linear interpolation are also shown (solid squares).
}
\end{figure}

Note that, in order to maintain the same user-defined error $\epsilon=10^{-2}$,
the number of terms $\nnodes$ increases with $\wmax$.
This relationship is exhibited more clearly in Figure~\ref{fig:Hvseps}, which shows the dependence between $\epsilon$ and $\nnodes$ for various $\wmax$.
This figure also shows the corresponding relationship for linear interpolation
(Section~\ref{sec:linear}).
For $\epsilon=10^{-2}$, linear interpolation requires \textapprox$10$ times as many nodes as \meth; this ratio is quite similar to the speedup shown in Figure~\ref{figure:results}.
Furthermore, we note that $H$ for the \meth method roughly follows the $\bigO((W + \log (1/\epsilon))^2)$ scaling given in Theorem~\ref{t:rank}.

\begin{figure}
\centering
  \includegraphics[width=.74\linewidth,trim= {0cm 0cm 0cm 0cm},clip]{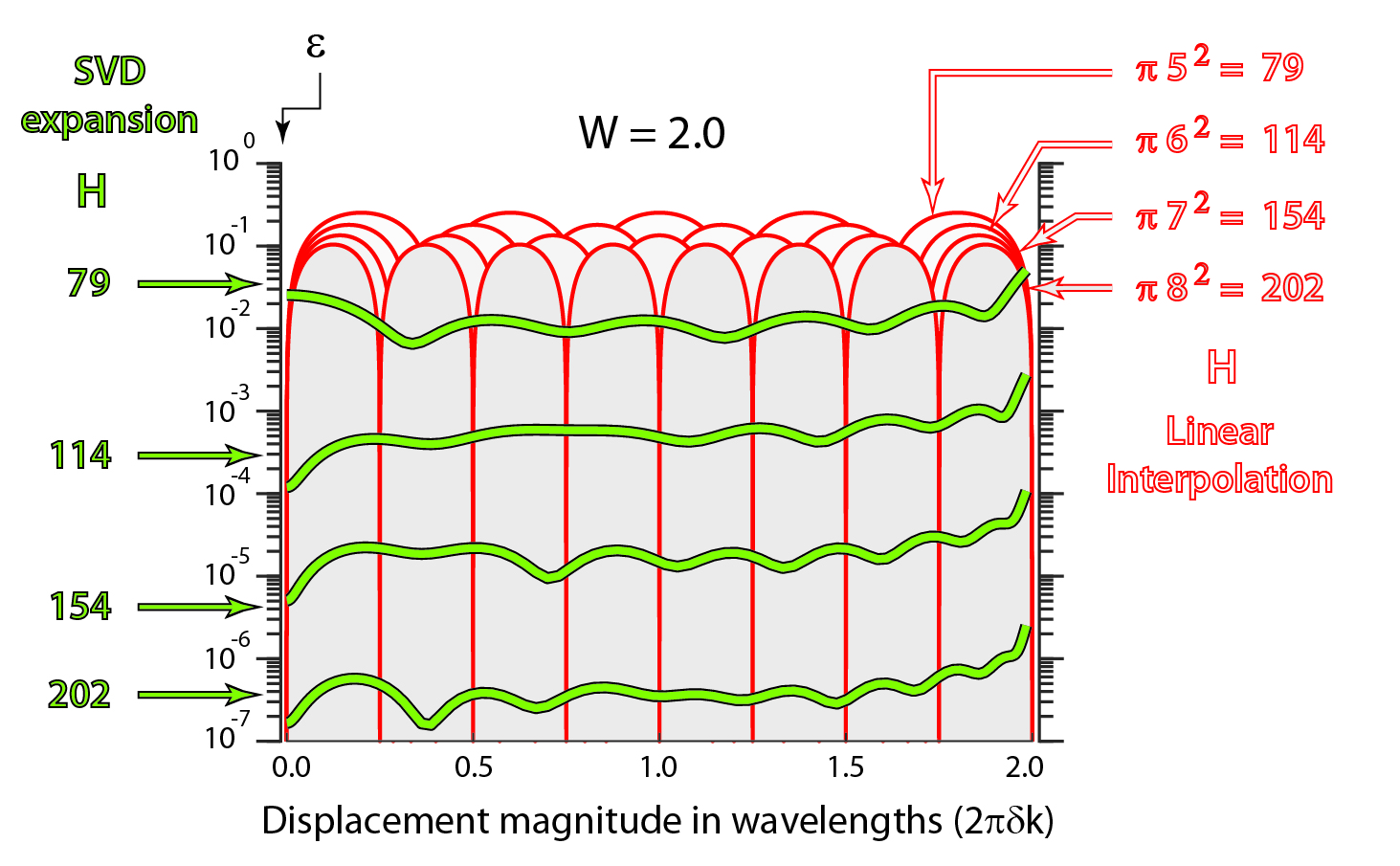}
\caption{ \label{fig:dkvseps}
The accuracy of our SVD expansion (solid green lines) as a function of the shift in wavelengths $2 \pi \delta K$ for $W=2$ using a semi-logarithmic scale.
We measure accuracy in terms of $\varepsilon(\delta)=\|F(\vd, \cdot) - F^\svd(\vd, \cdot)\|_2$ (i.e., the 2-norm on $\Omega_{\kmax}$), where the displacement direction $\omega$ is fixed at zero.
The value $\log_{10}(1/\varepsilon)$ can be interpreted as the number of digits of accuracy of our expansion.
Each solid green line corresponds to one particular choice of $\nnodes$ for our SVD expansion (values listed to the left).
Results for linear interpolation are also shown (red lines, number of nodes listed to the right).
}
\end{figure}

Finally, we emphasize that, as Table~\ref{table:complexity} makes clear, in order for \meth to be more efficient than BFT, the number of required terms $\nnodes$ must be smaller than the requested number of translations $\ntrans$, which requires a sufficient density of translations.
What is not explicitly stated in the table is that a coarsely sampled grid of translations may not be sufficient to reconstruct the landscape of inner products.
Indeed, while BFT produces accurate inner products for the translations on the grid, we must interpolate for off-grid translations.
As shown in Figure~\ref{fig:dkvseps}, a popular scheme such as linear interpolation does not accurately reproduce these intermediate inner products.
By contrast, \meth maintains nearly uniform error across the space of translations $\delta\leq\dmax$ for the accuracy determined by $\nnodes$, as explained in Section~\ref{sec:erroranalysis}.

\section{Extensions \label{sec:ext}}

\subsection{Other 2D kernels}

Taking a step back, we can describe our approach in broad terms.
Our original goal was to calculate inner products across a variety of translations.
The first step was to notice that the translation operator $T_\vd$ corresponds to pointwise multiplication by $F(\vd, \vk)$ over the Fourier domain $\vk$.
Our method boils down to approximating $F(\vd, \vk)$ with a sum of separable terms, each of the form $Y_\zeta(\vd) \cdot G_\zeta(\vk)$.

The same approach can also be applied to other convolutional operators.
For example, we might consider the inner product between one image and another that has been translated by some $\vd$, rotated by some $\gamma$, and convolved with a filter which depends on some parameter vector $\vect{\tau}$.
This is the case in cryo-EM, where images are filtered by some contrast transfer function (CTF), which is parametrized by a set of defocus parameters \cite{MindellGrigorieff2003}.
In addition to aligning the image with a CTF-filtered template, we may want to find the defocus parameters that give the highest correlation, which can yield improved reconstruction accuracy \cite{BartesaghiEtal2018}.

In the Fourier domain, translating by $\vd$ and filtering by a CTF with defocus parameters $\vect{\tau}$ corresponds to pointwise multiplication by a function dependent on $\vd$ and $\vect{\tau}$.
We can extend our approach to tackle this scenario by factorizing this function into terms of the form $Y_\zeta(\vd, \vect{\tau}) \cdot G_\zeta(\vk)$.
Plugging this factorization into \meth then yields a fast algorithm for computing inner products as a function of $\vd$, $\gamma$, and $\vect{\tau}$.

\subsection{3D volume-to-volume registration}

We have so far discussed the calculation of inner products between 2D images.
A closely related problem is the calculation of inner products between 3D volumes for the purposes of volume-to-volume alignment.
In this context, rotations can be described by their Euler angles $\alpha$, $\beta$, and $\gamma$, and the translations $\vd$ are in $\Real^3$.
At first glance, it appears as though our methodology should generalize naturally from 2D to 3D.
Unfortunately, there are several important differences which make such a generalization difficult.

To understand this impasse, first recall that our overall strategy was to (a) choose a basis in which rotations can be applied cheaply, (b) factorize the action of the translation operator in that basis using the SVD, and finally (c) use this factorization to calculate the inner products ${\cal X}(\vd, \gamma)$ over a set of translations $\vd$ and rotations $\gamma$.

The first problem arises in step (a).
In 2D, rotations commute, and as a result, there exists a basis in which rotations act diagonally: the Fourier--Bessel basis.
However, rotations do not commute in 3D, so there is no such basis for volumes.
A compromise is the spherical Fourier--Bessel basis, where the angular component is given by spherical harmonic functions.
In this basis, rotations are block diagonal, allowing us to compute multiple rotations quickly using 2D FFTs.

Step (b) proceeds in the same way for 3D, but step (c) incurs an additional cost.
Indeed, once we have factorized the translation operator in the spherical Fourier--Bessel basis, we multiply each right singular function $V_\zeta(\vk)$ by the coefficients of the volumes in the basis.
The resulting expression corresponds to one inner product in 2D, but not in 3D.
Thus, if we use the spherical Fourier--Bessel basis, applying each term in the SVD expansion requires significantly more work to calculate than a single inner product.

Another option is the 3D generalization of the least-squares interpolant described in Section~\ref{sec:lsq}.
The normal equations are similar to \eqref{eq:normal2d}, with the associated matrix
\begin{equation}
\label{eq:normal3dA}
\besselker^\ls(\vd^\prime, \vd) = \frac{\sin(\tilde{\delta}\kmax) - \tilde{\delta}\kmax \cos(\tilde{\delta}\kmax)}{(\tilde{\delta}\kmax)^{3}}~{,}
\end{equation}
which has a derivative
\begin{equation}
\label{eq:normal3dB}
\partial_{\tilde{\delta}\kmax} \besselker^\ls(\vd^\prime, \vd) = \frac{ ((\tilde{\delta}\kmax)^{2} - 3) \sin(\tilde{\delta}\kmax) + 3(\tilde{\delta}\kmax) \cos(\tilde{\delta}\kmax) }{(\tilde{\delta}\kmax)^{4}}~\text{.}
\end{equation}
As before, we may optimize the node locations using gradient descent, but we encounter the same instability issues at high accuracies.

Alternatively, if one of the volumes is known beforehand (e.g., when aligning a given volume against a fixed set of volumes), then we can use a functional SVD to precompute a factorized representation of 
\begin{equation}
X(\vd,\alpha,\beta,\gamma) = \iiint_{\Real^3} (T_{\vd} R_{\alpha\beta\gamma} \fA(\vk))^\ast \fB(\vk) \mathop{d\vk}
\end{equation}
such that:
\begin{equation}
X(\vd,\alpha,\beta,\gamma) \approx \sum_{\zeta=1}^\nnodes Y_{\zeta}(\vd) \iiint_{\Real^3} G_{\zeta}(\vk;\alpha,\beta,\gamma) \fB(\vk) \mathop{d\vk}
\end{equation}
where $R_{\alpha\beta\gamma}$ represents 3D rotation by the Euler angles $\alpha,\beta,\gamma$.
Here, $Y_1, \ldots, Y_\nnodes$ and $G_1, \ldots, G_\nnodes$ depend on the specifics of volume $A$, but not on $B$.

\section{Conclusion}
\label{sec:conc}

We introduced a new efficient method---named \meth---for calculating multiple inner products between two images of $n\times n$ pixels, one fixed and another rotated and translated with respect to the first.
The method is
based on decomposing the images in a Fourier--Bessel basis, allowing us to compute the inner products quickly for multiple rotation angles via 1D FFTs.
Furthermore, the truncated SVD of the translation operator enables us to compute approximate inner products efficiently for a large number $\ntrans$ of
arbitrary translations lying in a disk of radius $2W$ pixels.
The method has computational complexity $\bigO(H n (n+\ntrans))$,
where we prove (Theorem~\ref{t:rank})
that the rank $H = \bigO((W + \log(1/\epsilon))^2)$,
where $\epsilon$ is the relative root mean square error.
The method is evaluated on a set of synthetic cryo-EM projections
with sub-pixel translation grid spacings,
where it is shown to significantly outperform the competing BFR and BFT methods.
At mean spacing $\dx/4$ we show an acceleration over other known inner-product-forming methods of around $10\times$ over a wide range of $W$.
The method is $\bigO(n \log n)$ times faster than
the popular FFT-based cross-correlation alignment method (BFR), making it
favorable for large images.
Finally we present possible extensions of the method to other 2D kernels and 3D volume alignment, noting that the latter case poses additional difficulties due to the non-commutativity of rotations in 3D.

\ack

The authors thank
Alex Townsend for showing us the expansion \eqref{wimp} crucial to
Lemma~\ref{l:rank}, and Amit Singer, Leslie Greengard, and Zydrunas Gimbutas
for fruitful discussions and helpful suggestions.
The Flatiron Institute is a division of the Simons Foundation.

\section*{References}

\bibliographystyle{abbrv}  
\bibliography{frtalign_bib}

\end{document}